\documentclass[12pt]{article}
\usepackage{amssymb,amsfonts,amsmath,psfrag,eepic,colordvi,graphicx,epsfig}
\usepackage{amssymb,latexsym,graphics,array}
\usepackage{MnSymbol}
\usepackage[enableskew]{youngtab}
\usepackage{float,enumerate,enumitem}
\usepackage{tikz,ifthen}
\usepackage{pgfplots}
\usetikzlibrary{math}

\usepackage{framed}
\usepackage{hyperref}
\hypersetup{colorlinks=true}

\topskip  
\newcommand{\rmnum}[1]{\romannumeral #1}

\numberwithin{equation}{section}
\newtheorem{theorem}{Theorem}[section]

\newtheorem{conjecture}[theorem]{Conjecture}
\newtheorem{remark}[theorem]{Remark}
\newtheorem{lemma}[theorem]{Lemma}

	\def\des{\mathsf{des}}
\def\maj{\mathsf{maj}}
\def\exc{\mathsf{exc}}
\def\den{\mathsf{den}}
\def\st{\mathsf{st}}
\def\Inv{\mathsf{Inv}}
\def\Des{\mathsf{Des}}
\def\Exc{\mathsf{Exc}}
\def\Nexc{\mathsf{Nexc}}
\def\Nexcp{\mathsf{Nexcp}}
\def\Excp{\mathsf{Excp}}
\def\Excl{\mathsf{Excl}}
\def\inv{\mathsf{inv}}
\begin{document}
	\parskip 6pt
	
	\pagenumbering{arabic}
	\def\sof{\hfill\rule{2mm}{2mm}}
	\def\ls{\leq}
	\def\gs{\geq}
	\def\SS{\mathcal S}
	\def\qq{{\bold q}}
	\def\MM{\mathcal M}
	\def\TT{\mathcal T}
	\def\EE{\mathcal E}
	\def\lsp{\mbox{lsp}}
	\def\rsp{\mbox{rsp}}
	\def\pf{\noindent {\it Proof.} }
	\def\mp{\mbox{pyramid}}
	\def\mb{\mbox{block}}
	\def\mc{\mbox{cross}}
	\def\qed{\hfill \rule{4pt}{7pt}}
	\def\block{\hfill \rule{5pt}{5pt}}
	\def\lr#1{\multicolumn{1}{|@{\hspace{.6ex}}c@{\hspace{.6ex}}|}{\raisebox{-.3ex}{$#1$}}}
	\def\red{\textcolor{red}}

\begin{center}
	{\Large\bf Further results on $r$-Euler-Mahonian statistics}
\end{center}

\begin{center}
 
	{\small Kaimei Huang,  
		Sherry H.F. Yan$^{*}$\footnote{$^*$Corresponding author.}\footnote{{\em E-mail address:} hfy@zjnu.cn. }}

  Department of Mathematics,
 Zhejiang Normal University\\
 Jinhua 321004, P.R. China

\end{center}

\noindent {\bf Abstract.}
As natural generalizations of the descent number  ($\des$)   and the major index ($\maj$), Rawlings introduced the notions of the  $r$-descent number ($r\des$) and  the $r$-major index ($r\maj$)  for a given positive integer $r$.  A pair   $(\st_1, \st_2)$  of permutation statistics is said to be   $r$-Euler-Mahonian  if  $  (\mathrm{st_1}, \mathrm{st_2})$ and 
   $ (r\des,  r\maj)$  are equidistributed over the set  $\mathfrak{S}_{n}$ of all permutations of $\{1,2,\ldots, n\}$.  The main objective of this paper is to confirm   a recent conjecture  posed by Liu which asserts that     $(g\exc_\ell, g\den_\ell)$  is $(g+\ell-1)$-Euler-Mahonian for all positive integers $g$ and $\ell$, where $g\exc_\ell$ denotes the $g$-gap $\ell$-level excedance number and  $g\den_\ell$ denotes the $g$-gap $\ell$-level Denert's statistic. This is accomplished via a bijective proof of the equidistribution of $(g\exc_\ell, g\den_\ell)$ and $ (r\des,  r\maj)$ where $r=g+\ell-1$.    
    Setting $g=\ell=1$, our  result recovers the equidistribution of $(\des, \maj)$ and $(\exc, \den)$, which was first conjectured by Denert and proved by Foata and Zeilberger. 
Our second main result is concerned with the analogous result for  $(g\exc_\ell, g\den_{g+\ell})$ which states  that  $(g\exc_\ell, g\den_{g+\ell})$  is $(g+\ell-1)$-Euler-Mahonian for all positive integers $g$ and $\ell$.




\section{Introduction}
   Let $\mathfrak{S}_{n}$ denote the set of permutations of $[n]:=\{1,2,\ldots, n\}$. A permutation $\pi\in \mathfrak{S}_{n}$ is usually written in one-line notation as $\pi=\pi_1\pi_2\cdots \pi_n$. Given a permutation  $\pi\in \mathfrak{S}_{n}$, a pair $(i,j)$ is called an {\em inversion } of $\pi$ if $i<j$ and $\pi_i>\pi_j$. Let $\inv(\pi)$ denote the number of inversions of $\pi$.    
    Define the {\em descent set} of $\pi$ as 
    $$
    \Des(\pi)=\{i\mid \pi_i>\pi_{i+1},\,\, 1\leq i\leq n-1\}. 
   $$
   The {\em descent number } of $\pi$, denoted by $\des(\pi)$, is defined to be the cardinality of the set $\Des(\pi)$. 
  The {\em major index} of $\pi$ is defined as
$$
\maj(\pi)=\sum\limits_{i\in \Des(\pi) } i.
$$
The well-known MacMahon's equidistribution  theorem \cite{Mac} states that
\begin{equation}\label{invmaj}
	\sum\limits_{\pi\in \mathfrak{S}_{n}}  q^\mathrm{inv(\pi)}=\sum\limits_{\pi\in  \mathfrak{S}_{n }}q^\mathrm{maj(\pi)}=[n]_q[n-1]_q\ldots [1]_q,
\end{equation}
where $[k]_q:=1+q+q^2+\ldots+q^{k-1}$. 
Any permutation statistic equidistributed with $\inv$ is said to be  {\em Mahonian}. 
Various interesting  extensions of MacMahon's equidistribution theorem have been obtained   for other combinatorial objects  (see e.g. \cite{Chen1, Chen2, Liu1, Liu2, Wilson, Yan}).

Let $r\geq 1$. 
Rawlings~\cite{Raw} introduced the notions of {\em $r$-descent number} and   {\em $r$-major
index}, which are natural generalizations  of the  descent number and the major index. 
   Let $r\geq 1$. Given a permutation $\pi=\pi_1\pi_2\ldots \pi_n\in \mathfrak{S}_{n}$,  define the {\em $r$-descent set}  and the {\em $r$-inversion set} of $\pi$   to be
$$
r\Des(\pi)=\{i\in [n-1]\mid \pi_{i}\geq \pi_{i+1}+r\},
$$
$$
r\Inv(\pi)=\{(i,j)\in \Inv(\pi)\mid  \pi_{i}<\pi_{j}+r  \},
$$
respectively, where $\Inv(\pi)=\{(i,j)\mid \pi_i>\pi_j, \,\, i<j\}$.
   The {\em $r$-descent  number} of $\pi$, denoted by $r\des(\pi)$, is defined to be the   cardinality of the set $r\Des(\pi)$. 
The {\em   $r$-major
	index} of $\pi$ is defined to be 
$$
r\maj(\pi)=\sum\limits_{i\in r\Des(\pi)}i +|r\Inv(\pi)|. 
$$
  Clearly,  the statistic     $r\maj$ interpolates between  $\maj$ and $\inv$ since $r\maj$ reduces to $\maj$ when $r=1$, and reduces to $\inv$ when $r\geq n$. 
Rawlings~\cite{Raw}  presented a bijective proof of the equidistribution of  $r\maj$ and $\inv$ over $\mathfrak{S}_{n}$, thereby proving that $r\maj$ is Mahonian. A pair  $(\st_1, \st_2)$  of permutation statistics is said to be {\em $r$-Euler-Mahonian} if 
$$
\sum\limits_{\pi\in \mathfrak{S}_{n}}  t^{\st_1(\pi)}q^{\st_2(\pi)}=  \sum\limits_{\pi\in \mathfrak{S}_{n}}  t^{r\des(\pi)}q^{r\maj(\pi)}.
$$

 Given a permutation  $\pi=\pi_1\pi_2\ldots \pi_n\in \mathfrak{S}_{n}$,
an index $i$, $1\leq i< n$, is called an {\em excedance place}
 of $\pi$ if $\pi_i>i$. The {\em excedance number} of $\pi$, denoted by $\exc(\pi)$,  is defined to be the number of excedance places of $\pi$.    
 Denert \cite{Denert} introduced a new permutation statistic   $\den$ and conjectured that   $(\exc, \den) $  has the same distribution as  $(\des, \maj)$  over $\mathfrak{S}_{n}$, that is,
 \begin{equation}\label{den-maj}
 	 \sum_{\pi\in \mathfrak{S}_{n}}t^{\des(\pi)}q^{\maj(\pi)}=\sum_{\pi\in \mathfrak{S}_{n}}t^{\exc(\pi)}q^{\den(\pi)}.
 \end{equation}
  This conjecture was first  proved by Foata and Zeilberger \cite{Foata-Zeilberger}  and  a bijective proof was provided by Han \cite{Han}.

  Motivated by finding extensions of (\ref{den-maj}) for $r$-Euler-Mahonian statistics, 
    Liu~\cite{Liujcta} introduced the notions of {\em  $g$-gap $\ell$-level excedance number} and  {\em $g$-gap $\ell$-level Denert's statistic}. Fix   integers $g, \ell\geq 1$.  Given a permutation $\pi\in \mathfrak{S}_{n}$, an index $i\in[n-1]$ is called a  {\em $g$-gap $\ell$-level excedance place}  and $\pi_i$ is called a  {\em $g$-gap $\ell$-level excedance letter} if $\pi_i\geq i+g$ and $\pi_i\geq \ell$. For example, if we let    $\pi=836295417$,  then the $2$-gap  $3$-level excedance places of $\pi$ are given by $1,3, 5$ and the  $2$-gap  $3$-level excedance letters of $\pi$ are given by $6,8,9$. 
   Let $g\Excp_\ell(\pi)$ and $g\Nexcp_\ell(\pi)$ denote the set of $g$-gap $\ell$-level excedance places of $\pi$ and the set of non-$g$-gap-$\ell$-level excedance places of $\pi$, respectively.  Let $g\Excl_\ell(\pi)$ denote the set of $g$-gap $\ell$-level excedance letters of $\pi$. 
     Assume that  $g\Excp_\ell(\pi)=\{i_1, i_2, \ldots, i_k\}$ with $i_1<i_2<\cdots<i_k$ and   $g\mathrm{Nexcp}_\ell(\pi)=\{j_1, j_2, \ldots, j_{n-k}\}$
       with $j_1<j_2<\cdots<j_{n-k}$.  Define
   $ 
   g\Exc_{\ell}(\pi)=\pi_{i_1}\pi_{i_2}\ldots \pi_{i_k},
   $ 
   and 
   $ 
  g\Nexc_{\ell}(\pi)=\pi_{j_1}\pi_{j_2}\ldots \pi_{j_{n-k}}.
   $ 
    The {\em  $g$-gap $\ell$-level Denert's statistic} of $\pi$, denoted by $g\den_\ell(\pi)$, is defined as
   $$
   g\den_{\ell}(\pi)=\sum\limits_{i\in g\Excp_\ell(\pi)} (i+g-1)+{\inv}(g\Exc_{\ell}(\pi))+{\inv}(g\Nexc_{\ell}(\pi)).
   $$
    For example, if we let    $\pi=836295417$, then $2\Excp_3(\pi)=\{1,3,5\}$,   $2\Excl_3(\pi)=\{6,8,9\}$, $2\mathrm{Nexcp}_3(\pi)=\{2,4,6,7,8,9\}$, 
   $2\Exc_{3}(\pi) =869$, and  $2\Nexc_{3}(\pi) =325417$. Therefore,   we have
   $$
   \begin{array}{ll}
   	2\den_3(\pi)&=\sum\limits_{i\in 2\Excp_3(\pi)} (i+g-1)+{\inv}(2\Exc_{3}(\pi))+{\inv}(2\Nexc_{3}(\pi))\\
   	&= 2+4+6+{\inv}( 869)+{\inv}(325417)\\
      	&= 19.
   \end{array}
   $$   
   When $g=\ell=1$, the statistic $g\den_{\ell}$ reduces to the  Denert's statistic $\den$ introduced by Denert \cite{Denert}.
  The $g$-gap $1$-level Denert's statistic  of $\pi$ is called the {\em $g$-gap    Denert's statistic}  of $\pi$, denoted by $g\den(\pi)$. Similarly, The $1$-gap $\ell$-level Denert's statistic  of $\pi$ is called the {\em $\ell$-level     Denert's statistic } of $\pi$, denoted by $\den_\ell(\pi)$. 
     
  The {\em $g$-gap $\ell$-level excedance number } of $\pi$,  denoted by $g\exc_\ell(\pi)$,  is defined to be the number of $g$-gap $\ell$-level excedance letters that occurs weakly to the right of $\pi_{\ell}$, that is, 
    $$g\exc_\ell(\pi)=|\{i\mid \pi_i\geq i+g, i\geq \ell\}|.$$  For example, if we let    $\pi=836295417$,  then we  have
    $2\exc_{3}(\pi) =2$.  
  The $g$-gap $1$-level  excedance number of $\pi$ is called the {\em $g$-gap     excedance number}  of $\pi$, denoted by $g\exc(\pi)$. Similarly, the $1$-gap $\ell$-level  excedance number  of $\pi$ is called the {\em $\ell$-level      excedance number } of $\pi$, denoted by $\exc_\ell(\pi)$. 
  Recently, Liu   \cite{Liujcta} proved that the pair $(g\exc, g\den)$ is $g$-Euler-Mahonian for all $g\geq 1$ and further posed the following two conjectures.
  
  \begin{conjecture}[Liu~\cite{Liujcta}, Conjecture 1]\label{con1}     
  	For   $\ell\geq 1$, the pair $(\exc_\ell, \den_\ell)$ is $\ell$-Euler-Mahonian.
  \end{conjecture}
     
   \begin{conjecture}[Liu~\cite{Liujcta}, Conjecture 2]\label{con2}     
   	For all $g\geq 1$ and $\ell\geq 1$, the pair $(g\exc_\ell, g\den_\ell)$ is $r$-Euler-Mahonian, where $r=g+\ell-1$.
   \end{conjecture}
 It should be mentioned that Conjecture \ref{con1} has been confirmed by  Huang-Lin-Yan \cite{Huang-Lin-Yan}. 
 Our second main result    is concerned with the following analogous result of Conjecture \ref{con2}.
 \begin{theorem}\label{thm2}     
	For all $g\geq 1$ and $\ell\geq 1$, the pair $(g\exc_\ell, g\den_{g+\ell})$ is $r$-Euler-Mahonian, where $r=g+\ell-1$.
\end{theorem}

 From Table \ref{table1}, one can easily see that 
 $$
 \begin{array}{lll}
 	\sum\limits_{\pi\in \mathfrak{S}_{4}}  t^{3\des(\pi)}q^{3\maj(\pi)}&=& 1+3q+5q^2+5q^3+3q^4+q^5+tq^3+2tq^4+2tq^5+tq^6\\
 	&=&  \sum\limits_{\pi\in \mathfrak{S}_{4}}  t^{ 2\exc_{2}}q^{2\den_{4}}.
 \end{array}
 $$\begin{table*}\label{table1}
 	\begin{center}
 		\caption{The values of the pairs ($3\des$,$3\maj$) and ($2\exc_2$,$2\den_4$) of each permutation in $\mathfrak{S}_4$.}
 		\vspace{1em}
 		\fontsize{10}{8}\selectfont
 		\begin{tabular}{|c|c|c|}
 			\hline
 			\boldmath{$\mathfrak{S}_4$}& \boldmath{($3\des$,$3\maj$)}& \boldmath{($2\exc_2$,$2\den_4$)}\\
 			\hline
 			1234 &  (0,0) & (0,0) \\
 			1243 &  (0,1) & (0,1) \\
 			1324 &	(0,1) &	(0,1) \\
 			1342 &	(0,2) &	(0,2) \\
 			1423 &	(0,2) &	(1,3) \\
 			1432 &	(0,3) &	(1,4) \\
 			2134 &	(0,1) &	(0,1) \\
 			2143 &	(0,2) &	(0,2) \\
 			2314 &	(0,2) &	(0,2) \\
 			2341 &	(1,5) &	(0,3) \\
 			2413 &	(1,4) &	(1,4) \\
 			2431 &	(0,3) &	(1,5) \\
 			3124 &	(0,2) &	(0,2) \\
 			3142 &	(0,3) &	(0,3) \\
 			3214 &	(0,3) &	(0,3) \\
 			3241 &	(1,6) &	(0,4) \\
 			3412 &	(1,5) &	(1,5) \\
 			3421 &	(0,4) &	(1,6) \\
 			4123 &	(1,3) &	(0,2) \\
 			4132 &	(1,4) &	(0,3) \\
 			4213 &	(0,3) &	(0,3) \\
 			4231 &	(0,4) &	(0,4) \\
 			4312 &	(0,4) &	(0,4) \\
 			4321 &	(0,5) &	(0,5) \\
 			\hline
 		\end{tabular}
 		
 	\end{center}
 	
 \end{table*}

In this paper, we shall present bijective proofs of 
Conjecture \ref{con2} and Theorem \ref{thm2} by introducing 
two insertion methods   in the spirit of  Liu's bijective proof of the equidistribution of $(g\exc, g\den)$ and 
$(g\des, g\maj)$.

   \section{Rawlings's and Liu's insertion methods}     
     Let $r\geq 1$.
    In order to prove the equidistribution of $  r\maj$ and $\inv$,     
     Rawlings \cite{Raw} introduced an extension of 
      Carlitz's insertion method  \cite{Carlitz}. Here, we give a brief discription of this insertion  method. 
    
    Given $\pi\in \mathfrak{S}_{n-1}$, let 
   $$
   S(\pi)=\{\pi_{i}\mid \pi_{i-1}\geq \pi_{i}+r\}\cup \{\pi_{i}\mid \pi_i>n-r\}. 
   $$
  It is easily seen that  $|S(\pi)|=r\des(\pi)+r-1$. The {\em $r$-maj-labeling} of $\pi$  is obtained as follows.
   \begin{itemize}
   	\item Label the  space after $\pi_{n-1}$ by $0$.
   	\item Label the spaces before those letters of $\pi$ that are belonging to $S(\pi)$ from right to left with $1,2,\ldots, r\des(\pi)+r-1$. 
   	\item Label the remaining spaces from left to right with $r\des(\pi)+r, \ldots, n-1$.
   \end{itemize}

Take    $\pi=836295417$ and $r=3$ for example. Clearly,  we  have   $S(\pi)=\{1, 2, 3, 5, 8, 9  \}$,  and thus  the $3$-maj-labeling of $\pi$ is given by
$$
\begin{array} {llllllllll}
_{6}8&_53&_7 6&_4 2&_39&_2 5&_8 4&_11&_97&_0,
\end{array}
$$
where the labels of the spaces are written as subscripts. 
Define $$\phi^{\maj}_{r,n}:\mathfrak{S}_{n-1}\times  \{0, 1, \ldots, n-1\} \longrightarrow \mathfrak{S}_{n}$$ by mapping the pair $( \pi, c)$ to the  permutation in $\mathfrak{S}_{n}$  obtained by inserting   $ n$ at the space in $\pi$ which  is labeled by $c$ in the $r$-$\maj$-labeling of $\pi$. Continuing with our running example, we have
$\phi^{\maj}_{3,10}(836295417, 4)=836(10)295417$.

Rawlings \cite{Raw} proved that the map $\phi^{\maj}_{r,n}$  verifies  the following celebrated  properties.

\begin{lemma}[Rawlings~\cite{Raw}] \label{Raw}
	The map $\phi^{\maj}_{r,n}:\mathfrak{S}_{n-1}\times  \{0, 1, \ldots, n-1\} \longrightarrow \mathfrak{S}_{n}$ is a bijection  such that for $\pi\in \mathfrak{S}_{n-1}$ and $0\leq c\leq n-1$, we have 
	$$
	r\des(\phi^{\maj}_{r,n}(\pi, c))=\left\{ \begin{array}{ll}
		r\des(\pi)&\, \mathrm{if}\,\, 0\leq c\leq r\des(\pi)+r-1,\\
			r\des(\pi)+1&\, \mathrm{otherwise}, 
		\end{array}
	\right.
	$$
	and 
	$$
	r\maj(\phi^{\maj}_{r,n}(\pi, c))=r\maj(\pi)+c.
	$$
\end{lemma}

 Now we proceed to give a description of   Liu's  insertion method \cite{Liujcta}, which was employed to  prove that the pair   $(g\exc, g\den)$ is $g$-Euler-Mahonian for all $g\geq 1$. 

A $g$-gap $1$-level excedance letter is called a $g$-gap  excedance letter. 
 Given $\pi\in \mathfrak{S}_{n-1}$, assume that $g\exc(\pi)=s$ and that $e_1<e_2<\cdots<e_s$ are the $g$-gap excedance letters of $\pi$. The {\em $g$-gap-den-labeling of the index sequence} of $\pi$ is obtained by labelling the spaces of the index sequence $12 \cdots (n-1)$ of $\pi$ as follows.
 \begin{itemize}
	\item Label the  rightmost $g$ spaces  by $0,1,\ldots, g-1$.
	\item Label the space  before the index $e_i-(g-1) $ with $g+s-i$ for all $1\leq i\leq s$.
	\item Label the remaining spaces from left to right with $g+s, \ldots, n-1$.
\end{itemize}
The {\em  insertion letters} of the index sequence of $\pi$ are obtained by labeling  the space  to the left of the index $e_1-(g-1)$ with $e_1$, labeling  the space   to the left of the index $e_2-(g-1)$ and to the right of the index $e_1-(g-1)$ with $e_2$, $\ldots$, and labeling the space to the right of the index $e_s-(g-1)$ with $n$. 

For example,   we let    $\pi=836295417$.  Then the $3$-gap excedance letters are given by  $e_1=6$, $e_2=8$ and  $e_3=9$ and  thus the $3$-gap-den-labeling  and the  insertion letters    of the index sequence of  $\pi$ are given by$$
\begin{array}{lllllllllll} 
	&_{6}1&_{7}2&_{8}{ 3}&_{5}4&_{9}5&_{4}{\bold   6}&_{3}7&_{2}{\bold 8}&_{1}{\bold 9}&_{0}\\
	& \uparrow&\uparrow &\uparrow &\uparrow &\uparrow&\uparrow&\uparrow&\uparrow&\uparrow&\uparrow\\
	&6 &6 &6 &6 &8&8&9&10&10&10,
\end{array}
$$
where the labels of the spaces in the $3$-gap-den-labeling are written as subscripts and the insertion letters of the spaces  are written on the bottom row.

Define $$\phi^{\den}_{g,n}: \mathfrak{S}_{n-1}\times \{0, 1, \ldots, n-1\} \longrightarrow \mathfrak{S}_{n}$$ by mapping the pair $(  \pi, c)$ to the  permutation in $\mathfrak{S}_{n}$  according to the following rules.
\begin{itemize}
 \item If $c=0$, then let $\phi^{\den}_{g,n}(\pi,c)$ be the permutation obtained from $\pi$ by inserting $n$ immediately  after $\pi_{n-1}$. 
 \item For $c\geq 1$, assume that the space labeled by $c$ in the $g$-gap-den-labeling of the index sequence of $\pi$ is immediately before the index $x$, and assume that the insertion letter for this space is $e_{y}$, where $1\leq y\leq g\exc(\pi)+1$ and $e_{g\exc(\pi)+1}=n$. Define $\pi'=\pi'_1\pi'_2\ldots \pi'_{n-1}$ from $\pi$ by replacing $e_j$ with  $e_{j+1}$ for all $y\leq j\leq g\exc(\pi) $. Then define $\phi^{\den}_{g,n}(\pi,c)$ to be the permutation obtained from $\pi'$ by inserting the letter $e_y$ immediately before the letter $\pi'_{x}$.  
\end{itemize}

Take     $\pi=836295417$,  $g=3$ and $c=4$  for   example. Clearly, the space labeled by $4$ in the $3$-gap-den-labeling of the index sequence of $\pi$ is immediately before the index $x=6$, and   the insertion letter for this space is $e_2=8$.  Then we get a permutation $\pi'=9362(10)5417  $ from $\pi$ by replacing replacing $e_2=8$ with $e_3=9$ and  replacing $e_3=9$ with $e_4=10$. Finally, we get $\phi^{\den}_{3,10}(\pi,4)= 9362(10)85417$  from $\pi'$ by inserting $e_2=8$ immediately before $\pi'_x=\pi'_6=5$.

Liu~\cite{Liujcta} proved that the map $\phi^{\den}_{g,n}$ possesses  the following desired properties.

\begin{lemma}[Liu~\cite{Liujcta}] \label{Liu}	 
	The map $\phi^{\den}_{g,n}: \mathfrak{S}_{n-1}\times \{0, 1, \ldots, n-1\} \longrightarrow \mathfrak{S}_{n}$ is a bijection  such that 
	for $\pi\in \mathfrak{S}_{n-1}$ and $0\leq c\leq n-1$, we have
		$$ 
		g\exc(\phi^{\den}_{g,n}(\pi))=\left\{ \begin{array}{ll}
			g\exc(\pi)&\, \mathrm{if}\,\, 0\leq c\leq g\exc(\pi)+g-1,\\
			g\exc(\pi)+1&\, \mathrm{otherwise}, 
		\end{array}
		\right.
		$$
		and 
		$$
		g\den(\phi^{\den}_{g,n}(\pi))=g\den(\pi)+c.
		$$

\end{lemma}
Relying on Lemmas \ref{Raw} and  \ref{Liu},  Liu \cite{Han} provided a bijective proof of the equidistribution of  $(g\exc, g\den)$ and $(g\des, g\maj)$.

\section{Proof of Conjecture \ref{con2}}

This section is devoted to the proof of Conjecture \ref{con2}. Our proof is based on introducing a  nontrivial extension of Liu's insertion method \cite{Liujcta}. 

Let
  $n\geq g+\ell$ and $g,\ell\geq 1$. 
Given $\pi\in \mathfrak{S}_{n-1}$, assume that  $g\Excl_\ell(\pi)=\{e_1, e_2, \ldots, e_s\}$ with $e_1<e_2<\cdots<e_s$.
Define 
$$
\mathcal{A}(\pi)=\{i\mid    \pi_i\in g\Excl_\ell(\pi), \,\,   \pi_i\geq g+\ell,\,\, i<\ell\}, 
$$
$$
\mathcal{B}(\pi)=\{i\mid  \pi_i\notin g\Excl_\ell(\pi), \,\,    i<\ell\}, 
$$
and 
$$
\mathcal{C}(\pi)=\{i\mid \pi_i\in g\Excl_\ell(\pi), \,\,   \pi_i<g+\ell,\,\, i<\ell\}.
$$
It is apparent that
$|\mathcal{A}(\pi)|+|\mathcal{C}(\pi)|+g\exc_\ell(\pi)=s$ and $|\mathcal{A}(\pi)|+|\mathcal{B}(\pi)|+|\mathcal{C}(\pi)|=\ell-1$. Thus, we have $|\mathcal{B}(\pi)|+s=g\exc_\ell(\pi)+\ell-1$.   

 The {\em $g$-gap-$\ell$-level-den-labeling} of    the index sequence  of $\pi$ is obtained by labeling the spaces of the index sequence $12\ldots (n-1)$ of $\pi$ as follows.
\begin{itemize}
	\item Label the  rightmost $g$ spaces  by $0,1,\ldots, g-1$ from right to left.
\item   For all $1\leq i\leq s$,  label the space  before the index $e_i-(g-1) $ with $g+s-i$ if $e_i\geq g+\ell$, and label the space  before the index $x $ with $g+s-i$ if   $e_i<g+\ell$ and $\pi_x=e_i$. 
		\item 	Label the spaces  before the indices   belonging to $\mathcal{B}(\pi)$ from left to  right with $g+s, \ldots, g\exc_\ell(\pi)+g+\ell-2$.
	\item Label the remaining spaces from left to right with $g\exc_\ell(\pi)+g+\ell-1, \ldots, n-1$.
\end{itemize}

Assume that the set   $g\Excl_{g+\ell}(\pi)$ of $g$-gap $(g+\ell)$-level excedance letters of $\pi$ is given by  $$g\Excl_{g+\ell}(\pi)=\{e_i\mid e_i\geq g+\ell, 1\leq i\leq s\}=\{e_t, e_{t+1}, \ldots, e_s\}.$$ If $g\Excl_{g+\ell}(\pi) \neq \emptyset$,  
the {\em $g$-gap $(g+\ell)$-level insertion letters}  of the index sequence of $\pi$ are obtained by labeling  the spaces  to the left of the index $e_t-(g-1)$ and   to the right of the index $\ell-1$ with $e_t$, labeling  the spaces   to the left of the index $e_{t+1}-(g-1)$ and to the right of the index $e_{t}-(g-1)$ with $e_{t+1}$, $\ldots$, and labeling the spaces to the right of the index $e_s-(g-1)$ with $n$.  
If $g\Excl_{g+\ell}(\pi) =\emptyset$,   the {\em  $g$-gap $(g+\ell)$-level insertion letters} of the index sequence of $\pi$ are obtained by labeling   the spaces to the right of the index $\ell-1$ with $n$.

For example, let    $\pi= 596284317$, $g=3$ and $\ell=6$.  Clearly, we have $\mathcal{A}(\pi)=\{2\}$,  $\mathcal{B}(\pi)=\{1,4\}$,  $\mathcal{C}(\pi)=\{3,5\}$, $3\Excl_{6}(\pi)=\{6,8,9\}$ and  $3\Excl_{9}(\pi)=\{9\}$. Then the $3$-gap-$6$-level-den-labeling and the $3$-gap $9$-level insertion letters    of  the index sequence of  $\pi$ are  given by$$
\begin{array}{lllllllllll} 
	&_{6}1&_{8}2&_{5}3&_{7}4&_{4}5&_{9}{\bold 6}&_{3}7&_{2}{\bold 8}&_{1}{\bold 9}&_{0}\\
		&  &  &  &  & &\uparrow &\uparrow&\uparrow&\uparrow&\uparrow\\
	&  &  &  &  & &9&9&10&10&10,
		 
\end{array}
$$
where the labels of the spaces   are written as subscripts and the insertion letters of the spaces  are written on the bottom row.

\begin{theorem}\label{thden}
Let $g, \ell\geq 1$ and $n\geq g+\ell$.
 There is a bijection
  $$\phi^{\den}_{n,g,\ell}:  \mathfrak{S}_{n-1}\times \{0, 1, \ldots, n-1\} \longrightarrow \mathfrak{S}_{n}$$ such that  for $\pi\in \mathfrak{S}_{n-1}$ and $0\leq c\leq n-1$,  we have 
  
  \begin{equation}\label{eqprop1}
  g\exc_\ell(\phi^{\den}_{n,g,\ell}( \pi, c))=\left\{ \begin{array}{ll}
  	g\exc_\ell(\pi)&\, \mathrm{if}\,\, 0\leq c\leq g\exc_\ell(\pi)+g+\ell-2,\\
  	g\exc_\ell(\pi)+1&\, \mathrm{otherwise}, 
  \end{array}
  \right.
  \end{equation}
  and 
  \begin{equation}\label{eqprop2}
  g\den_\ell(\phi^{\den}_{n,g,\ell}( \pi, c))=g\den_\ell(\pi)+c.
  \end{equation} 
\end{theorem}

In the following, we shall introduce  four maps   which will play
essential roles in the construction of $\phi^{\den}_{n,g,\ell}$.

For any positive integers $a$ and $b$ with $a< b$,  let $(a,b)$ denote the set of all integers $i$ with $a<i<b$,  and let $[a,b)$ denote the set of all integers $i$ with $a\leq i<b$.
For $n\geq \ell$, 
an index $i\in[n]$ is called a {\em $g$-gap  $\ell$-level  grande-fixed place} of $\tau=\tau_1\tau_2\ldots \tau_n\in\mathfrak{S}_n$ if one of the following conditions hold.
\begin{itemize}
	\item $\tau_i=n$ with $n-(g-1)< i\leq n$;
	\item    $ \tau_i=i+(g-1) $ and   $\tau_i\geq \ell$;
	
	\item      $\tau_i\in  g\Excl_{\ell}(\tau)$ and  $[i+g-1, \tau_i)\cap g\Excl_{\ell}(\tau)=\emptyset$.
\end{itemize}
 For example, if  $\tau=836279415$, $g=3$ and $\ell=6$, then we have $3\Excl_{6}(\tau)=\{6, 8, 9\}$. It is easily checked  that the indices $3$ and $5$ are $3$-gap $6$-level grande-fixed places, while the indices $1$ and $6$ are not  $3$-gap $6$-level grande-fixed places.

 Let $\mathfrak{S}_{n, g, \ell}$ denote the subset of permutations   $\tau\in \mathfrak{S}_{n}$  such that there does not exist any $g$-gap $(g+\ell)$-level  grande-fixed place $z$ of $\tau$  with $z\geq \ell$. 
   Denote by 
  $ 
  \mathfrak{S}^{(1)}_{n,g, \ell}
  $  the subset of  permutations $\tau\in \mathfrak{S}_{n, g, \ell}$  with  $\tau_{\ell}\geq g+\ell$.

  Let $\tau=\tau_1\tau_2\ldots \tau_n $ be a permutation in $\mathfrak{S}_{n, g, \ell}-\mathfrak{S}^{(1)}_{n,g, \ell}$.   Choose  $\tau_{x}$  to be the greatest $g$-gap $\ell$-level excedance letter  located   to the left of $\tau_\ell$.  Find all the non-$g$-gap-$\ell$-level excedance letters of $\tau$  located between $\tau_{x}$ and $\tau_{\ell}$ (including $\tau_{\ell}$), say $\tau_{i_0}, \tau_{i_1}, \ldots, \tau_{i_k}$ with $\ell=i_0>i_1>i_2>\cdots >i_k>i_{k+1}=x$.  Then  $\tau$ is said to be a {\em proper}  permutation in  $\mathfrak{S}_{n, g, \ell}-\mathfrak{S}^{(1)}_{n,g, \ell}$ if   we have either    $\tau_{i_{j}}<\ell$  or  $  \tau_{i_{j}}<i_{j+1}+g$  for all $0\leq j\leq k$.
Otherwise, it is said to be an {\em improper} permutation in   $\mathfrak{S}_{n, g, \ell}-\mathfrak{S}^{(1)}_{n,g, \ell}$.  For example, the permutation $ \tau=(10)975864312$ is a proper permutation in  $\mathfrak{S}_{10, 3, 6}-\mathfrak{S}^{(1)}_{10, 3, 6}$, whereas the permutation $\tau'=(10)965874312$ is an improper permutation in   $\mathfrak{S}_{10, 3, 6}-\mathfrak{S}^{(1)}_{10, 3, 6}$.
 Let $\mathfrak{S}^{(2)}_{n, g, \ell}$ (resp. $\mathfrak{S}^{(3)}_{n, g, \ell}$)   denote  the set of  proper (resp. improper)  permutations in $ \mathfrak{S}_{n, g, \ell}-\mathfrak{S}^{(1)}_{n, g, \ell}$.

Let $\mathcal{X}_{n,g, \ell}$ (resp., $\mathcal{Y}_{n,g, \ell}$ and  $\mathcal{Z}_{n,g, \ell}$) denote  the subset of pairs $(\pi, c)\in \mathfrak{S}_{n-1}\times \{ 1, \ldots,  n-1\}$ such that the space labeled by $c$  under the $g$-gap-$\ell$-level-den-labeling of the index sequence of $\pi$    is immediately  before the index $x$   for some $ x\in\mathcal{A}(\pi)$ (resp.,  $ x\in\mathcal{B}(\pi)$ and  $ x\in\mathcal{C}(\pi)$ ).

 \begin{framed}
 \begin{center}
{\bf The map $\alpha_{n,g, \ell}:\mathcal{X}_{n,g, \ell}\longrightarrow \mathfrak{S}^{(1)}_{n, g, \ell}$} 
\end{center}
 Let $(\pi, c)\in \mathcal{X}_{n,g, \ell}$.  Suppose   that
 $$ g\Excl_{\ell}(\pi) =\{e_1, e_2, \ldots, e_s\}$$
 and 
   $$  g\Excl_{g+\ell}(\pi)=\{e_i\mid e_i\geq g+\ell, 1\leq i\leq s\}=\{e_t, e_{t+1}, \ldots, e_s\}$$ with $e_1<e_2<\cdots<e_s$.  Assume that   the space labeled by $c$  under the $g$-gap-$\ell$-level-den-labeling of the index sequence of $\pi$    is  immediately before the index $x$.   
Suppose that $\mathcal{A}(\pi)=\{i_1, i_2, \ldots, i_k\}$ where $i_1<i_2<\cdots<i_k$ and  $x=i_{y}$ for some $1\leq y\leq k$. 
Define  $\alpha_{n, g, \ell}(\pi,c)$ to be the permutation  constructed by the following procedure.
\begin{itemize}
	\item Step 1: Construct a permutation  $\sigma=\sigma_1\sigma_2\ldots \sigma_{n-1}$  from $\pi$ by replacing $e_j$ with  $e_{j+1}$ for all $t\leq j\leq s$ with the convention that $e_{s+1}=n$. 
	\item Step 2:  Define   $\alpha_{n, g, \ell}(\pi,c)$ to be the permutation  obtained  from $\sigma$ by replacing $\sigma_{i_y}=\sigma_x$ with $e_t$,
	replacing $\sigma_{i_j}$ with $\sigma_{i_{j-1}}$ for all $y<j\leq k$,   and inserting $ \sigma_{i_k}$ immediately before $\sigma_{\ell}$.
\end{itemize}
 \end{framed}
Continuing with our running example where $\pi= 596284317$, $g=3$ and $\ell=6$,   we have $e_1=6$, $e_2=8$,   $e_3=9$, $3\Excl_9(\pi)=\{9\}$ and    $\mathcal{A}(\pi)=\{2\}$.
Recall that the $3$-gap-$6$-level-den-labeling    of the index sequence of  $\pi$ is given by$$
\begin{array}{lllllllllll} 
	&_{6}1&_{8}2&_{5}3&_{7}4&_{4}5&_{9}{\bold 6}&_{3}7&_{2}{\bold 8}&_{1}{\bold 9}&_{0}. 
	
\end{array}
$$
 It is easily seen that the space labeled by $8$ in the    $3$-gap-$6$-level-den-labeling  of the index sequence of $\pi$ is immediately before the index $x=2 \in \mathcal{A}(\pi)$. First we generate a permutation $\sigma=5(10)6284317$ from $\pi$ by replacing $e_3=9$ with $e_4=10$.  Then we get the permutation $\alpha_{10, 3,6}(\pi,8)=\tau= 59628(10)4317 $   by replacing $\sigma_2=10$ with $e_3=9$ and      inserting $\sigma_2=10$ immediately before $\sigma_{6}=4$.

Now we proceed to   show that the map $\alpha_{n,g, \ell}$ verifies the following desired properties. 
\begin{lemma}\label{lemalpha1}
	Let $g, \ell\geq 1$ and $n\geq g+\ell$.  For any
	$(\pi, c)\in \mathcal{X}_{n,g, \ell} $, we have 
		\begin{equation}\label{alphare1}
		g\exc_\ell(\alpha_{n,g, \ell}(\pi, c))= g\exc_\ell(\pi)+1 	 	 
		\end{equation}
	and	 
			\begin{equation}\label{alphare2}
		g\den_{g+\ell}(\alpha_{n,g, \ell}(\pi, c))-g\den_{g+\ell}(\pi)=g\den_\ell(\alpha_{n,g, \ell}(\pi, c))-g\den_\ell(\pi)=c.
		 \end{equation}	 
 \end{lemma}
\pf Here we use the  notations of the definition of  the map $\alpha_{n,g, \ell}$. Let $\alpha_{n,g, \ell}(\pi, c)=\tau=\tau_1\tau_2\ldots \tau_n$.
Recall that $\sigma$ is obtained from $\pi$ by replacing $e_i$ with $e_{i+1}$ for all $t\leq i\leq s$ with the convention $e_{s+1}=n$. Since   $e_i<e_{i+1}$ for all $t\leq i\leq s+1$,    it follows that the procedure from $\pi$ to $\sigma$ preserves  the    set of $g$-gap $(g+\ell)$-level excedance places and  leaves all the non-$g$-gap-$(g+\ell)$-level excedance letters in their place. Moreover, this procedure keeps the relative order of   the $g$-gap $(g+\ell)$-level excedance letters as well as the relative order of   the $g$-gap $\ell$-level excedance letters.   Therefore, we have $ 
g\den_\ell(\sigma)=g\den_\ell(\pi) 
$, $g\exc_{\ell}(\pi)=g\exc_{\ell}(\sigma)$ and  $ 
 g\den_{g+\ell}(\sigma)=g\den_{g+\ell}(\pi) 
 $.

 In the following, we shall verify that the following facts hold.\\
 \noindent{\bf Fact 1:} The set of   $g$-gap  $(g+\ell)$-level excedance letters located   to the left of $\sigma_{\ell}$
 is  given by $\{\sigma_{i_j}\mid  1\leq j\leq k\}$. \\
  \noindent{\bf Fact 2:}  The set of the  $g$-gap  $(g+\ell)$-level excedance letters located weakly to the left of $\tau_{\ell}$  is  given by $\{\tau_{i_j}\mid  1\leq j\leq k+1\}$ with the convention that $i_{k+1}=\ell$.\\
  \noindent{\bf Fact 3:}   The letter $\sigma_i$ is a $g$-gap  $\ell$-level excedance letter (resp. $g$-gap  $(g+\ell)$-level excedance letter ) of $\sigma$ if and only if $\tau_i$ is a $g$-gap  $\ell$-level excedance letter (resp. $g$-gap  $(g+\ell)$-level excedance letter )  of $\tau$ for all $i<\ell$.  Moreover, the letter $\tau_{\ell}$ is a   $g$-gap  $(g+\ell)$-level excedance letter   of $\tau$. \\
  \noindent{\bf Fact 4:} For all $i\geq \ell$, the letter $\sigma_i$ is a $g$-gap  $\ell$-level excedance letter  of $\sigma$ if and only if $\tau_{i+1}=\sigma_i$ is a $g$-gap  $\ell$-level excedance letter  of $\tau$.   
   
 Recall that  $$\mathcal{A}(\pi)=\{i\mid   \pi_i\geq g+\ell, \,\,i<\ell\}=\{i_1, i_2, \ldots, i_k\}.$$ 
Then Fact 1 follows immediately from the fact   that the procedure from $\pi$ to $\sigma$ preserves  the    set of $g$-gap $(g+\ell)$-level excedance places. 
 Recall that $\tau$ is obtained from $\sigma$ by replacing $\sigma_x=\sigma_{i_y}$ with $e_t$,
 replacing $\sigma_{i_j}$ with $\sigma_{i_{j-1}}$ for all $y<j\leq k$,   and inserting $ \sigma_{i_k}$ immediately before $\sigma_{\ell}$.  Since $e_t\geq g+\ell$  and   $\sigma_{i_j}\geq g+\ell$  for all $1\leq j\leq k$ by Fact 1, 
 it follows that the letter  $\tau_{i_j}$ is also a $g$-gap $(g+\ell)$-level excedance letter of $\tau$ for all $1\leq j\leq k+1$, completing the proof of  Fact 2. 
 By the construction of $\tau$, we have
 $\tau_z=\sigma_{z}$ if $z<\ell$ and  $z\notin \{i_j\mid 1\leq j\leq k\}$. Then Fact 3 follows directly from Fact 1 and Fact 2. 

By the construction of $\sigma$,  one can see that   $\sigma_i\in g\Excl_{g+\ell}(\sigma)$ would imply that $\sigma_i>\pi_i\geq i+g$. 
   By the construction of $\tau$,  we have $\tau_{i+1}=\sigma_i$ for all $i\geq \ell$. 
Hence, we deduce that
$\sigma_i\geq i+g$   if and only if $\tau_{i+1}\geq i+1+g$  for all $i\geq \ell$,   completing the proof of Fact 4.

  Fact  3  and  Fact 4  tell us   that the procedure from $\sigma$ to $\tau$ will increase $g\exc_\ell(\sigma)$ by $1$. As $g\exc_\ell(\pi)=g\exc_\ell(\sigma)$, we have
 $ 
 g\exc_\ell(\tau)=g\exc_\ell(\pi)+1,
 $ 
 completing the  proof of~\eqref{alphare1}.
 
 Next we proceed to prove ~\eqref{alphare2}.
Recall that 
$$
g\den_\ell(\sigma)=\sum\limits_{i\in g\Excp_\ell(\sigma)} (i+g-1)+{\inv}(g\Exc_{\ell}(\sigma))+{\inv}(g\Nexc_{\ell}(\sigma)) 
$$ 
and 
$$
g\den_{g+\ell}(\sigma)=\sum\limits_{i\in g\Excp_{g+\ell}(\sigma)} (i+g-1)+{\inv}(g\Exc_{g+\ell}(\sigma))+{\inv}(g\Nexc_{g+\ell}(\sigma)). 
$$     
 By Fact 3 and  Fact 4,   one can easily check that the procedure from $\sigma$ to $\tau$   increases   the sum $\sum\limits_{i\in g\Excp_\ell(\sigma)}(i+g-1)$ (resp.  the sum $\sum\limits_{i\in g\Excp_{g+\ell}(\sigma)}(i+g-1)$) by 
 $\ell+g-1+g\exc_\ell(\sigma)$. 
 
By Fact 1 through Fact 4,  one can easily see  that the procedure from $\sigma$ to $\tau$  keeps the relative order of  the non-$g$-gap-$\ell$-level excedance letters  as well as  the relative order of   the non-$g$-gap-$(g+\ell)$-level excedance letters. Therefore, this  procedure  does not affect   $\inv(g\Nexc_\ell(\sigma))$ and $\inv(g\Nexc_{g+\ell}(\sigma))$. 
    
    By the construction of $\sigma$, one can see that  each  $g$-gap $(g+\ell)$-level excedance letter of $\sigma$ is greater than $e_t$.  
  Fact 1 tells us  that the number  $g$-gap $(g+\ell)$-level excedance letters of $\sigma$  located to the left of $\sigma_{i_y}$ is given  by $y-1$. Therefore, the replacement of  $\sigma_{i_y}$ with $e_t$ in the procedure from $\sigma$ to $\tau$ will increase   $\inv(g\Exc_\ell(\sigma))$  (resp. $\inv(g\Exc_{g+\ell}(\sigma))$)  by $y-1$. 
      Note that $g\den_\ell(\sigma)=g\den_\ell(\pi)$, $g\den_{g+\ell}(\sigma)=g\den_{g+\ell}(\pi)$, and $g\exc_\ell(\sigma)=g\exc_\ell(\pi)$.
    Hence, we derive  that
    $$
    	g\den_\ell(\tau)-g\den_\ell(\pi)=g\den_\ell(\tau)-g\den_\ell(\sigma)=g\exc_\ell(\sigma)+g+\ell+y-2= g\exc_\ell(\pi)+g+\ell+y-2 
    $$
    and 
   $$
   \begin{array}{ll}
   g\den_{g+\ell}(\tau)-g\den_{g+\ell}(\pi)&=g\den_{g+\ell}(\tau)-g\den_{g+\ell}(\sigma)\\
   &=g\exc_\ell(\sigma)+g+\ell+y-2\\
   &= g\exc_\ell(\pi)+g+\ell+y-2. 
   \end{array}
   $$ 
 According to the rules specified in  the $g$-gap-$\ell$-level-den-labeling    of the index sequence of  $\pi$, we have $c=  g\exc_\ell(\pi)+g+\ell-2+y$. 
 Therefore, we conclude that 	$g\den_{g+\ell}(\tau)-	g\den_{g+\ell}(\pi)=g\den_\ell(\tau)-g\den_\ell(\pi)=c$, completing the  proof of~\eqref{alphare2}.
    \qed

 \begin{lemma}\label{lemalpha2}
 	Let $g, \ell\geq 1$ and $n\geq g+\ell$.  The map  $\alpha_{n, g,\ell}  $ is an injection from $\mathcal{X}_{n,g,\ell}$ to $\mathfrak{S}^{(1)}_{n, g, \ell}$. 
 \end{lemma}
 \pf Here we use the  notations of the definition of  the map $\alpha_{n,g, \ell}$. 
  Let $\alpha_{n,g, \ell}(\pi, c)=\tau=\tau_1\tau_2\ldots \tau_n$. 
First, we aim to   show  that $\tau\in \mathfrak{S}^{(1)}_{n, g, \ell} $.  
 By Fact~3, we have $\ell\in g\Excp_{g+\ell}(\tau)  $ and thus $\tau_{\ell}\geq g+\ell$. Hence, in order to prove that $\tau\in \mathfrak{S}^{(1)}_{n,g, \ell} $, it remains to   show that  there does not exist any  $g$-gap   $(g+\ell)$-level   grande-fixed place $z$ of $\tau$ with $z\geq \ell$.  To this end, we need to verify the following fact.\\ 
  \noindent{\bf Fact 5:} $g\Excl_{g+\ell}( \tau)=\{e_i\mid t\leq i\leq s+1\}$ with the convention that $e_{s+1}=n$.\\
   By the  the construction of  $\sigma$, it is easily checked  that $g\Excl_{g+\ell}( \sigma)=\{ e_i\mid t<i\leq s+1\}. $
  By Fact 1 through  Fact 4, the procedure from $\sigma$ to $\tau$   transforms  $e_t$ into a $g$-gap $(g+\ell)$-level excedance letter   and preserves all the $g$-gap $(g+\ell)$-level excedance letters of $\sigma$.  Therefore,  we deduce that $g\Excl_{g+\ell}( \tau)=\{e_i\mid t\leq i\leq s+1\}$ as desired, completing the proof of Fact~5. 
  
   Fact~5 tells us that $e_{s+1}=n$ is a 
 $g$-gap  $(g+\ell)$-level excedance letter of $\tau$.  This implies that $\tau_i\neq n$ for all $n-(g-1)< i\leq n$.
 Hence,  the index $z$ is not a $g$-gap   $(g+\ell)$-level   grande-fixed place of $\tau$ for all $z>n-(g-1)$.

 Take $\tau_z=z+(g-1)$ arbitrarily (if any) with $z>\ell$.  By Fact~4, $\tau_z= z+(g-1)$ would imply  that $\sigma_{z-1}\leq z-1+(g-1)$. On the other hand, by the definition of the procedure from $\sigma$ to $\tau$,  we have  $ \sigma_{z-1}=\tau_z=z+(g-1)$,  a contradiction. Therefore,  we have $\tau_z\neq z+(g-1)$ for all $z> \ell$.

 Take $\tau_z\geq z+g$ arbitrarily (if any) with $z>\ell$.
 Again by Fact~4, the letter $\tau_{z}=\sigma_{z-1}$ is a $g$-gap  $(g+\ell)$-level excedance letter in $\sigma$.    By the construction of $\sigma$,  the letter $\pi_{z-1}$ is also a $g$-gap  $(g+\ell)$-level excedance letter in $\pi$.
 Again by the construction of $\sigma$, we have $\sigma_{z-1}>\pi_{z-1}\geq z-1+g$ and $\pi_{z-1} \in \{e_t, e_{t+1}, \ldots, e_s\}$. This implies that
 $ 
\pi_{z-1}\in  [z-1+g, \sigma_{z-1})\cap \{e_t, e_{t+1}, \ldots, e_s\}\neq \emptyset 
 $.
 Recall  that  $\tau_z=\sigma_{z-1}$ when $z>\ell$  and  $g\Excl_{g+\ell}( \tau)=\{e_t, e_{t+1}, \ldots, e_{s+1}\}$. 
 Then   we have   $[z+g-1, \tau_{z})\cap g\Excl_\ell(\tau)=[z-1+g, \sigma_{z-1})\cap g\Excl_\ell(\tau)\neq \emptyset $.  Hence, we conclude that  
 the index $z$ is not  a $g$-gap $(g+\ell)$-level   grande-fixed place of  $\tau$ for any $z> \ell$.

      By the construction of $\tau$, we have $\tau_{\ell}=\sigma_{i_k}$.  By Fact 1, the letter $\sigma_{i_k}$ is a $g$-gap $(g+\ell)$-level excedance letter of $\sigma$. By the construction of $\sigma$, we have $\sigma_{i_k}>e_t\geq g+\ell$. Then, by Fact 5, it follows that   $e_t\in[\ell+g-1, \sigma_{i_k}) \cap g\Excl_{g+\ell}(\tau)= [\ell+g-1, \tau_\ell) \cap g\Excl_{g+\ell}(\tau)\neq \emptyset $. This implies that the index $\ell$ is not a $g$-gap $(g+\ell)$-level grande-fixed place  of $\tau$. 
   So far, we have concluded that  
  the index $z$ is not  a $g$-gap $(g+\ell)$-level   grande-fixed place of  $\tau$ for any $z\geq  \ell$ and thus $\tau\in \mathfrak{S}^{(1)}_{n, g, \ell}$. 
 
 Next, we proceed to show that the map  $\alpha_{n, g,\ell}  $ is an injection.
  By (\ref{alphare2}), we have $c=g\den_\ell(\tau)-g\den_\ell(\pi)$. Therefore,  in order to  show that $\alpha_{n, g, \ell}$ is injective, it suffices to show that we can recover $\pi$ from $\tau$.  First we aim to show  the resulting permutations $\sigma$ and $\tau$ verifies the following desired properties.
  
\begin{itemize}
		\item[\upshape(\rmnum{1})] $g\Excl_{g+\ell}(\sigma)=\{e_i\mid t<i\leq  s+1 \}$ with the convention that $e_{s+1}=n$.
	\item[\upshape(\rmnum{2})]   $\tau_x=e_t$ is the smallest $g$-gap  $(g+\ell)$-level excedance letter of $\tau$.
	\item[\upshape(\rmnum{3})]  The set of the $g$-gap  $(g+\ell)$-level excedance letters that occur  between $\tau_x$ and    $\tau_\ell$ (including $\tau_x$ and $\tau_{\ell}$)  is given by 
	$\{\tau_{i_j}\mid y\leq j\leq k+1\}$  where  $x=i_{y}<i_{y+1}<\cdots< i_k<i_{k+1}=\ell$  and $\tau_{i_j}=\sigma_{i_{j-1}}$ for all $y<j\leq k+1$. 
 
\end{itemize} 
Property (\rmnum{1}) follows directly from the construction of $\sigma$ and 
  properties (\rmnum{2})-(\rmnum{3})  follow immediately   from   Fact~2 and Fact~5.

   Now we are in the position to give  a description of the procedure  to recover the permutation $\pi$ from $\tau$.
 \begin{itemize}
 	\item  Choose $\tau_{x'}$ to be  the smallest $g$-gap  $(g+\ell)$-level excedance letter of $\tau$. 
 	\item 
 	Find all the $g$-gap  $(g+\ell)$-level excedance letters,   say   $\tau_{j_1}, \tau_{j_2}, \ldots, \tau_{j_{k'}}$ with $x'=j_1<j_2<\cdots<j_{k'}=\ell$,  that are  located  between  $\tau_{x'}$ and   $\tau_\ell$  (including $\tau_{x'}$ and  $\tau_\ell$) in $\tau$.  
 	Construct a permutation $\sigma$ from $\tau$ by removing $\tau_{j_{k'}}$ and then replacing  each $\tau_{j_{i-1}}$  with   $\tau_{j_{i}}$ for all $1< i\leq k'$.
 	\item  Find all the $g$-gap  $(g+\ell)$-level  excedance  letters of $\sigma$, say $e'_1, e'_2, \ldots, e'_m$ with $e'_1<e'_2<\cdots<e'_m$.
   Let $\pi$ be the permutation obtained from $\sigma$ by replacing each $e'_i$ with $e'_{i-1}$ for all $1\leq i\leq m$ with the convention that $e'_0=\tau_{x'}$.
 \end{itemize} 
 Properties (\rmnum{1})-(\rmnum{3}) ensure  that after applying  the above procedure to $\tau$, we  will  recover the permutation  $\sigma$ and thus the   permutation $\pi$. This completes the proof. \qed

\begin{framed}
\begin{center}
{\bf The map $\beta_{n,g, \ell}:\mathcal{Y}_{n,g, \ell}\longrightarrow \mathfrak{S}^{(2)}_{n, g, \ell} $ }
\end{center}
  Let $(\pi, c)\in \mathcal{Y}_{n,g, \ell}$.  Suppose   that
 $$ g\Excl_{\ell}(\pi) =\{e_1, e_2, \ldots, e_s\}$$
 and 
 $$  g\Excl_{g+\ell}(\pi)=\{e_i\mid e_i\geq g+\ell,\,\, 1\leq i\leq s\}=\{e_t, e_{t+1}, \ldots, e_s\}$$ with $e_1<e_2<\cdots<e_s$.  Assume that   the space labeled by $c$  under the $g$-gap-$\ell$-level-den-labeling of the index sequence of $\pi$    is  immediately before the index $x$.   
 Suppose that $\mathcal{B}(\pi)=\{i_1, i_2, \ldots, i_k\}$ where $i_1<i_2<\cdots<i_k$ and  $x=i_{y}$ for some $1\leq y\leq k$. 
 Define  $\beta_{n,g,\ell}(\pi,c)$ to be the permutation  constructed by the following procedure.
 \begin{itemize}
 	\item Step 1: Construct a permutation  $\sigma=\sigma_1\sigma_2\ldots \sigma_{n-1}$  from $\pi$ by replacing $e_j$ with  $e_{j+1}$ for all $t\leq j\leq s$ with the convention that $e_{s+1}=n$. 
 	\item Step 2: Construct a permutation  
 	 $\theta=\theta_1\theta_2\ldots \theta_{n}$   from $\sigma$ by    replacing $\sigma_x=\sigma_{i_y}$ with $e_t$,
 	replacing $\sigma_{i_j}$ with $\sigma_{i_{j-1}}$ for all $y<j\leq k$,   and inserting $ \sigma_{i_k}$ immediately before $\sigma_{\ell}$.	
 
 	\item  Step 3: Find all the $g$-gap $(g+\ell)$-level excedance letters of $\theta$ that occur to the left of $\theta_\ell$, say $e'_1, e'_2, \ldots, e'_m$ with  $e'_1<e'_2<\cdots< e'_m$.  
 	Then define   $ \beta_{n, g, \ell}(\pi,c)$ to be the permutation obtained from $\theta$ by replacing
 	$e'_j$ with $e'_{j-1}$ for all $1< j\leq m$ and replacing $e'_1$ with $e'_m$. 
 \end{itemize}
 
  \end{framed}

 Continuing with our running example where $\pi= 596284317$, $g=3$ and $\ell=6$,   we have $e_1=6$, $e_2=8$,   $e_3=9$, $3\Excl_9(\pi)=\{9\}$  and $\mathcal{B}(\pi)=\{1,4\}$.
 Recall that the $3$-gap-$6$-level-den-labeling    of the index sequence of  $\pi$ is given by$$
 \begin{array}{lllllllllll} 
 	&_{6}1&_{8}2&_{5}3&_{7}4&_{4}5&_{9}{\bold 6}&_{3}7&_{2}{\bold 8}&_{1}{\bold 9}&_{0}. 
 	
 \end{array}
 $$  It is easily seen that the space labeled by $6$ in the   $3$-gap-$6$-level-den-labeling of the index sequence of $\pi$ is immediately before the index $x=1$ and $x\in \mathcal{B}(\pi)$.  First we generate a permutation $\sigma=5(10)6284317$ from $\pi$ by replacing $e_3=9$ with $e_4=10$.   Then we generate a permutation $\theta=9(10)65824317 $ from $\sigma$ by replacing $\sigma_1=5$ with $e_3=9$,  replacing $\sigma_4=2$ with $\sigma_1=5$,  and inserting $\sigma_4=2$ immediately before $\sigma_6=4$.  
 Clearly, all the $3$-gap $9$-level excedance letters that occur to the left of $\theta_6=4$ in $\theta$ are given by $e'_1=9$ and $ e'_2=10$.
 Finally, we get the permutation  $\beta_{10, 3,6}(\pi,6)=   (10)965824317$   obtained from $\theta$ by   replacing $e'_2=10$ with $e'_1=9$ and   replacing $e'_1=9$ with $e'_2=10$.

\begin{lemma}\label{lembeta1}
	Let $g, \ell\geq 1$ and $n\geq g+\ell$.  For any
$(\pi, c)\in \mathcal{Y}_{n,g, \ell} $, we have   
\begin{equation}\label{betare1}
	g\exc_\ell(\beta_{n,g, \ell}(\pi, c))= g\exc_\ell(\pi)		 	 
\end{equation}
and 
\begin{equation}\label{betare2}
	g\den_\ell(\beta_{n,g, \ell}(\pi, c))=g\den_\ell(\pi)+c.
\end{equation}
\end{lemma}
 
 \pf  Here we use the  notations of the definition of  the map $\beta_{n,g, \ell}$. Let $\beta_{n,g, \ell}(\pi, c)=\tau=\tau_1\tau_2\ldots \tau_n$. By the same reasoning as in the proof of Lemma \ref{lemalpha1}, one can verify that   $ 
 g\den_\ell(\sigma)=g\den_\ell(\pi) 
 $ and  $ 
 g\exc_\ell(\sigma)=g\exc_\ell(\pi) 
 $.

 In the following, we shall prove that the following facts hold. \\
 \noindent{\bf Fact ${\bf1'}$:} The set of   non-$g$-gap-$\ell$-level excedance letters located   to the left of $\sigma_{\ell}$
 is  given by $\{\sigma_{i_j}\mid  1\leq j\leq k\}$.\\
  \noindent{\bf Fact ${\bf2'}$:}  The set of the  non-$g$-gap-$\ell$-level excedance letters located weakly to the left of $\theta_{\ell}$  is  given by $\{\theta_{i_j}\mid  1\leq j\leq k+1,\,\, j\neq y\}$ with the convention that $i_{k+1}=\ell$.\\
  \noindent{\bf Fact ${\bf3'}$.} For all $i<\ell$ and $i\neq x$, the letter $\sigma_i$ is a $g$-gap  $\ell$-level excedance letter  of $\sigma$ if and only if $\theta_{i}$ is a $g$-gap  $\ell$-level excedance letter of $\theta$. Moreover,  the letter $\theta_x$ is a $g$-gap  $\ell$-level excedance letter of $\theta$, whereas  the letter  $\theta_\ell$ is  a non-$g$-gap-$\ell$-level excedance  letter of $\theta$.\\ 
 \noindent{\bf Fact ${\bf4'}$.} For all $i\geq \ell$, the letter $\sigma_i$ is a $g$-gap  $\ell$-level excedance letter of $\sigma$ if and only if $\theta_{i+1}=\sigma_i$ is a $g$-gap  $\ell$-level excedance letter  of $\theta$.\\      
  \noindent{\bf Fact ${\bf5'}$.} $g\Excl_{g+\ell}(\theta)=\{e_i\mid t\leq i\leq s+1\}$ with the convention that $e_{s+1}=n$.\\
 Fact $4'$ could be justified by  similar arguments as in the proof of  Fact~$4$ and the proof is omitted here. 
 Recall that  $$\mathcal{B}(\pi)=\{i\mid  \pi_i\notin g\Excl_{\ell}(\pi), \,\,  i<\ell\}=\{i_1, i_2, \ldots, i_k\}.$$ 
 Then Fact $1'$ follows immediately from the fact that  the procedure from $\pi$ to $\sigma$   preserves the set of $g$-gap $\ell$-level excedance places  and   leaves  all the non-$g$-gap- $\ell$-level excedance letters  in their place.  
   
 Recall that $\theta$ is obtained from $\sigma$ by replacing $\sigma_x=\sigma_{i_y}$ with $e_t$,
 replacing $\sigma_{i_j}$ with $\sigma_{i_{j-1}}$ for all $y<j\leq k$,   and inserting $ \sigma_{i_k}$ immediately before $\sigma_{\ell}$.  Since $e_t\geq g+\ell$, the letter $\theta_{i_y}=e_t$ is a $g$-gap $\ell$-level excedance letter of $\theta$.   By  Fact $1'$,    we have either $\sigma_{i_j}<\ell$ or $\sigma_{i_j}\leq i_j+g$ for all $1\leq j\leq k$. This implies that  
 the letter  $\theta_{i_j}$ is also    a non-$g$-gap-$\ell$-level excedance letter of $\theta$ for all $1\leq j\leq k+1$ and  $j\neq y$ by the construction of $\theta$, completing the proof of  Fact $2'$.    Since $\theta_z=\sigma_z$ if  $z<\ell$ and $z\notin \{i_j\mid 1\leq j\leq k\}$,  Fact $3'$   follows  directly from Fact $1'$ and Fact $2'$. 
  
By the construction of  $\sigma$, it is easily checked  that $g\Excl_{g+\ell}( \sigma)=\{e_i\mid t<i\leq s+1 \} $.
 Then, by Fact $1'$ through Fact $4'$,  one can easily check that the procedure from $\sigma$ to $\theta$   transforms  $e_t$ into a $g$-gap $(g+\ell)$-level excedance letter   and preserves all the $g$-gap $(g+\ell)$-level excedance letters of $\sigma$.   Therefore,  we deduce that $ g\Excl_\ell( \theta)=\{e_t, e_{t+1}, \ldots, e_{s+1}\}$ 
  as desired, completing the proof of Fact~$5'$.
So far, we have completed the proof   of Fact $1'$ through Fact $5'$.

 Clearly, the procedure from $\theta$ to $\tau$ preserves the number of   $g$-gap $\ell$-level excedance letters that occur weakly to the right of $\theta_\ell$ as $\tau_i=\theta_i$ for all $i\geq \ell$. By Fact $3'$ and  Fact $4'$, the procedure from $\sigma$ to $\theta$ also preserves the number of  $g$-gap $\ell$-level excedance letters that occur weakly to the right of $\sigma_\ell$. Hence, we derive that 
 $$
 g\exc_\ell(\tau) = g\exc_\ell(\theta)=g\exc_\ell(\sigma)= g\exc_\ell(\pi), 
 $$
 where the last equality follows from  the fact that $g\exc_\ell(\sigma)=g\exc_\ell(\pi)$, 
 completing the proof of (\ref{betare1}).

 In the following, we proceed to verify (\ref{betare2}).
 Recall that 
 $$
 g\den_\ell(\sigma)=\sum\limits_{i\in g\Excp_\ell(\sigma)} (i+g-1)+{\inv}(g\Exc_{\ell}(\sigma))+{\inv}(g\Nexc_{\ell}(\sigma)). 
 $$
  From  Fact $1'$ through Fact $4'$,   one can easily check that  the procedure from $\sigma$ to $\theta$  increases the sum $\sum\limits_{i\in \mathrm{gExcp}_\ell(\sigma)}(i+g-1)$ by $x+g-1+g\exc_\ell(\sigma)$ and  does not affect $\inv(g\Nexc_\ell(\sigma))$.

   Assume that there are  exactly $p$ $g$-gap $\ell$-level excedance letters and $p'$  $g$-gap $(g+\ell)$-level excedance letters   located  to the left of $\pi_x$,  and there are exactly  $q$ $g$-gap $\ell$-level excedance letters and      $q'$  $g$-gap $(g+\ell)$-level excedance letters  located between $\pi_x$ and $\pi_{\ell}$ in $\pi$.   
        By the rules specified  in the $g$-gap-$\ell$-level-den-labeling of the places of the index sequence of $\pi$, we have $c=(g-1)+s+y=(g-1)+g\exc_\ell(\pi)+p+q+y$ and $x=p+y$.
        Clearly, the procedure from $\pi$ to $\sigma$  preserves  the set of $g$-gap $ \ell$-level excedance places as well as   the set of $g$-gap $(g+\ell)$-level excedance places.  
           So, we conclude that there are  exactly $p$ $g$-gap $\ell$-level excedance letters and $p'$  $g$-gap $(g+\ell)$-level excedance letters   located  to the left of $\sigma_x$,  and there are exactly  $q$ $g$-gap $\ell$-level excedance letters and      $q'$  $g$-gap $(g+\ell)$-level excedance letters  located between $\sigma_x$ and $\sigma_{\ell}$ in $\sigma$.  
          Note that each $g$-gap $(g+\ell)$-level excedance letter of $\sigma$ is greater than $e_t$.  Then, the replacement of $\sigma_{x}=\sigma_{i_y}$ with $e_t$ in the procedure from $\sigma$ to $\theta$ will increase ${\inv}(g\Exc_{\ell}(\sigma))$ by
         $ p'+(q-q')$.        
           By the above analysis, we derive  that
   \begin{equation}\label{eqbeta1}
   	\begin{array}{ll}
   		g\den_\ell(\theta)&=\sum\limits_{i\in g\Excp_\ell(\theta)} (i+g-1)+{\inv}(g\Exc_{\ell}(\theta))+{\inv}(g\Nexc_{\ell}(\theta))\\
   		&=g\den_\ell(\sigma)+x+(g-1)+g\exc_\ell(\sigma)+p'+(q-q')\\
   		    		&=	g\den_\ell(\pi)+x+(g-1)+g\exc_\ell(\pi)+p'+(q-q'),\\
   	\end{array}
   \end{equation}
   where the last equality follows from the fact that $g\exc_\ell(\pi)=g\exc_\ell(\sigma)$ and $g\den_\ell(\pi)=g\den_\ell(\sigma)$.

 As $e'_i\geq g+\ell$ for all $1\leq i\leq m$,   the procedure from $\theta$ to $\tau$ preserves the set of $g$-gap $\ell$-level excedance  places as well as the set of $g$-gap $\ell$-level excedance  letters. Moreover, this procedure leaves all the non-$g$-gap-$\ell$-level excedance letters in their place.  Therefore,  this procedure   does not affect   the sum $\sum\limits_{i\in g\mathrm{
 		Excp}_\ell(\theta)}(i+g-1)$ and $\inv(g\Nexc_{\ell}(\theta))$.

 Recall that
 there are  exactly  $p'$  $g$-gap $(g+\ell)$-level excedance letters   located  to the left of $\sigma_x$,  and there are exactly      $q'$  $g$-gap $(g+\ell)$-level excedance letters  located between $\sigma_x$ and $\sigma_{\ell}$ in $\sigma$. Since the procedure from $\sigma$ to $\theta$ preserves the number of such $g$-gap $(g+\ell)$-level excedance letters, it follows that  
   the number of  $g$-gap $(g+\ell)$-level excedance letters   located  to the left of $\theta_x$  in $\theta$ is given by $p'$ and the number of 
 $g$-gap $(g+\ell)$-level excedance letters    located   between $\theta_x$ and   $\theta_\ell$ in $\theta$ is given by $q'$.    Fact $5'$ tells us that $\theta_x=e_t$ is the smallest $g$-gap $(g+\ell)$-level excedance letter of $\theta$. This yields that $\theta_x=e'_1$. 
 Then replacing $\theta_x=e'_1$ with $e'_m$ in the procedure from $\theta$ to $\tau$  will decrease $\inv(g\Exc_{\ell}(\theta)) $ by $p'$  
 and increase  $\inv(g\Exc_{\ell}(\theta)) $ by $q'$. Thus, we derive that 
 \begin{equation}\label{eqbeta3}
 	\inv(g\Exc_{\ell}(\tau))=\inv(g\Exc_{\ell}(\theta))-p'+q'.
 \end{equation}
 
  Recall that $g\exc_\ell(\sigma)=g\exc_\ell(\pi)$, $g\den_\ell(\sigma)=g\den_\ell(\pi)$, $x=p+y$, and   $c=(g-1)+g\exc_\ell(\pi)  +p+y+q$.
 Invoking (\ref{eqbeta1}) and (\ref{eqbeta3}), we derive that
 $$
 \begin{array}{ll}
 	g\den_\ell(\tau)&=\sum\limits_{i\in g\Excp_\ell(\tau)} (i+g-1)+{\inv}(g\Exc_{\ell}(\tau))+{\inv}(g\Nexc_{\ell}(\tau))\\
 	&=	g\den_\ell(\theta) -p'+q'\\ 
 	  	  	   	  	  		&= g\den_\ell(\pi)+(g-1) +g\exc_\ell(\pi)+x+q\\
 	  	  		&= g\den_\ell(\pi)+(g-1) +g\exc_\ell(\pi)+y+p+q\\
 	&= g\den_r(\pi)+c,
 \end{array}
 $$
 completing the proof of  (\ref{betare2}).   
  \qed

  \begin{lemma}\label{lembeta2}
  	Let $g, \ell\geq 1$ and $n\geq g+\ell$.  The map  $\beta_{n, g,\ell}  $ is an injection from $\mathcal{Y}_{n,g,\ell}$ to $\mathfrak{S}^{(2)}_{n, g, \ell}$. 
  \end{lemma}
 \pf   Here we use the  notations of the definition of  the map $\beta_{n,g, \ell}$.  Let $\beta_{n,g, \ell}(\pi, c)=\tau=\tau_1\tau_2\ldots \tau_n$.  First, we shall prove that  the resulting permutations $\tau$, $\theta$ and $\sigma$ possess the following desired properties.      
 \begin{itemize}
 		\item[\upshape(\rmnum{1})] $g\Excl_{g+\ell}(\sigma)=\{e_i\mid t<i\leq  s+1 \}$ with the convention that $e_{s+1}=n$.
 	\item[\upshape(\rmnum{2})]   $\theta_x=e_t$ is the smallest  $g$-gap $(g+\ell)$-level excedance letter of $\theta$. 
 		\item[\upshape(\rmnum{3})]  The set of the non-$g$-gap-$\ell$-level excedance letters   located   between  $\theta_x$ and   $\theta_\ell$ (including $\theta_\ell$) is given by 
 	$\{\theta_{i_j}\mid y<j\leq k+1\}$ where  $x<i_{y+1}<\cdots< i_{k+1}=\ell$ and $\theta_{i_j}=\sigma_{i_{j-1}}$ for all $y<j\leq k+1$.   		
   		\item[\upshape(\rmnum{4})]  The set of the  $g$-gap $(g+\ell)$-level excedance letters located    to the left of   $\tau_\ell$  is given by $\{e'_i\mid 1\leq i\leq m\}$. 
 	\item[\upshape(\rmnum{5})]   $\tau_x$ is the greatest  $g$-gap $\ell$-level excedance letter located  to the left of $\tau_\ell$ in $\tau$.	
 	 
 	\item[\upshape(\rmnum{6})]  The set of the non-$g$-gap-$\ell$-level excedance letters located   between  $\tau_x$ and   $\tau_\ell$ (including $\tau_\ell$) is given by 
 	$\{\tau_{i_j}\mid y<j\leq k+1\}$ where  $x<i_{y+1}<\cdots< i_{k+1}=\ell$ and $\tau_{i_j}=\sigma_{i_{j-1}}$ for all $ y<j\leq k+1$.  	 
 \end{itemize}  
Property (\rmnum{1}) follows directly from the construction of $\sigma$. Recall that $\theta_x=e_t$.  Then (\rmnum{2}) follows directly from Fact $5'$.   Property(\rmnum{3}) follows directly from Fact $2'$.  
   Since $e'_i\geq g+\ell$ for all $1\leq i\leq m$,  the procedure from $\theta$ to $\tau$ preserves  the set of $g$-gap $(g+\ell)$-level excedance letters located  to the left of $\theta_\ell$.  Hence, the set of  the $g$-gap $(g+\ell)$-level excedance letters located to the left of $\tau_{\ell}$ is given by $\{e'_i\mid 1\leq i\leq m\}$, completing the proof of (\rmnum{4}).  By (\rmnum{2}),  we have $e'_1=e_t=\theta_x$. Recall that in the procedure from $\theta$ to $\tau$, we replace $e'_1$ with $e'_m$. This yields that $\tau_x=e'_m$, completing the proof of  (\rmnum{5}).       Note that  the procedure from $\theta$ to $\tau$  only permutes the $g$-gap $(g+\ell)$-level excedance letters and leaves all the remaining  letters in their place. 
  Thus,       (\rmnum{6}) follows immediately from    (\rmnum{3}). So far, we have completed the proof of properties (\rmnum{1})-(\rmnum{6}).

 Next  we proceed to show that $\tau\in  \mathfrak{S}^{(2)}_{n, g, \ell} $.   By properties (\rmnum{5}) and (\rmnum{6}), in order to show $\tau$ is proper,  it suffices to show that  we have   either
 $\tau_{i_{j}}<\ell$ or 
 $   \tau_{i_{j}}<i_{j-1}+g$   for  all $y< j\leq k+1$.    By property (\rmnum{6}), we have $\tau_{i_j}=\sigma_{i_{j-1}}$ for all $y<j\leq k+1$.   Fact $1'$ tells us that    the letter $\sigma_{i_j}$ is a non-$g$-gap-$\ell$-level excedance letter of $\sigma$ for all $y\leq j\leq k$.  
 This implies that we have
 either $\sigma_{i_{j}}<\ell$ or 
 $  \sigma_{i_{j}}<i_{j}+g$  for  all $y\leq j\leq k$.   Hence, we deduce that we have either
 $\tau_{i_{j}}=\sigma_{i_{j-1}}<\ell$ or 
 $  \tau_{i_{j}}=\sigma_{i_{j-1}}<i_{j-1}+g$   for  all $y< j\leq k+1$.  So,  we have concluded that $\tau$ is indeed proper.

 In order to show that $\tau\in  \mathfrak{S}^{(2)}_{n, g, \ell} $, it remains to show that  $\tau\in  \mathfrak{S}_{n, g, \ell}-\mathfrak{S}^{(1)}_{n, g, \ell} $. 
 Clearly, we have $\tau_{\ell}=\theta_{\ell}$.
 By Fact $3'$, the index $\ell$ is not a $g$-gap $\ell$-level excedance place of $\theta$. Hence, we have $\tau_{\ell}< g+\ell$ and thus $\tau\notin \mathfrak{S}^{(1)}_{n, g, \ell}$. In order to show that $\tau\in  \mathfrak{S}_{n, g, \ell}-\mathfrak{S}^{(1)}_{n, g, \ell} $, it remains to show that $\tau\in \mathfrak{S}_{n, g, \ell}$.  Recall that $\tau_i=\theta_i$ for all $i\geq \ell$. In order to show that $\tau\in \mathfrak{S}_{n, g, \ell}$, it suffices to show  that    the index $i$ is not a $g$-gap  $(g+\ell)$-level grande-fixed place of $\theta$ for any $i\geq\ell$.   By the same reasoning as in the proof of    Lemma~\ref{lemalpha2}, one can verify that the index $i$ is not a $g$-gap  $(g+\ell)$-level grande-fixed place of $\theta$ for any $i>\ell$ relying on Fact $5'$.
     Recall that $\theta_{\ell}=\sigma_{i_k}$ by the construction of $\theta$. By Fact $1'$, we have $\sigma_{i_k}<g+\ell$. Hence, the index $\ell$ is not a
       $g$-gap  $(g+\ell)$-level grande-fixed place of $\theta$.  
 So far, we have  concluded that $\tau\in  \mathfrak{S}^{(2)}_{n, g, \ell} $.

 In the following, we aim to show that the map $\beta_{n, g, \ell}$ is injective. By (\ref{betare2}), we have $c=g\den_\ell(\tau)-g\den_\ell(\pi)$. In order to  show that $\beta_{n,g,\ell}$ is injective, it is sufficient  to show that we can recover $\pi$ from $\tau$.   Now we proceed to give  a description of the procedure  to recover the permutation $\pi$ from $\tau$.
 \begin{itemize}
  	\item  Find all the  $g$-gap  $(g+\ell)$-level excedance  letters of $\tau$  located  to the left of $\tau_\ell$, say $f_1, f_2, \ldots, f_{m'} $ with $f_1<f_2<\cdots<f_{m'}=\tau_{x'}$. 	Construct a permutation $\theta$ by replacing $f_{i-1}$ with $f_{i}$ for all $1<i\leq m'$ and replacing $f_{m'}$ with $f_1$. 
 	\item  
 	Find all the non-$g$-gap-$\ell$-level excedance letters,   say   $\theta_{j_1}, \theta_{j_2}, \ldots, \theta_{j_{k'}}$ with $x'<j_1<j_2<\cdots<j_{k'}=\ell$,  that are  located  to the right of $\theta_{x'}$ and weakly to the left of  $\theta_\ell$ in $\theta$.  
 	Construct a permutation $\sigma$ by removing $\theta_{j_{k'}}$ and then replacing  each $\theta_{j_{i-1}}$  with   $\theta_{j_{i}}$ for all $1\leq i\leq k'$ with the convention that $j_0=x'$. 
 	\item  Find all the  $g$-gap  $(g+\ell)$-level letters of $\sigma$, say $g_1, g_2, \ldots, g_h$ with $g_1<g_2<\cdots<g_h$.
 	Let $\pi$ be the permutation obtained from $\sigma$ by replacing  $g_i$ with $g_{i-1}$ for each $1\leq i\leq h$ with the convention that $g_0=\theta_{x'}$.
 \end{itemize} 
 Properties (\rmnum{1} )-(\rmnum{5} ) ensure that after applying  the above procedure to $\tau$, we  will recover  the permutation $\theta$,   the permutation  $\sigma$ and thus the permutation $\pi$. This completes the proof.\qed

 \begin{framed}
 	\begin{center}
 		{\bf The map $\gamma_{n,g, \ell}:\mathcal{Z}_{n,g, \ell}\longrightarrow \mathfrak{S}^{(3)}_{n, g, \ell} $ }
 	\end{center}
 	Let $(\pi, c)\in \mathcal{Z}_{n,g, \ell}$.  Suppose   that
 	$$ g\Excl_{\ell}(\pi) =\{e_1, e_2, \ldots, e_s\}$$
 	and 
 	$$  g\Excl_{g+\ell}(\pi)= \{e_i\mid e_i\geq g+\ell, 1\leq i\leq s\}=\{e_t, e_{t+1}, \ldots, e_s\}$$ with $e_1<e_2<\cdots<e_s$.  Assume that   the space labeled by $c$  under the $g$-gap-$\ell$-level-den-labeling of the index sequence of $\pi$    is  immediately before the index $x$ for some  $x\in \mathcal{C}(\pi)$.   
  	Define  $\gamma_{n,g,\ell}(\pi,c)$ to be the permutation  constructed by the following procedure.
 	\begin{itemize}
 		\item Step 1: Construct a permutation  $\sigma=\sigma_1\sigma_2\ldots \sigma_{n-1}$  from $\pi$ by replacing $e_j$ with  $e_{j+1}$ for all $t\leq j\leq s$ with the convention that $e_{s+1}=n$.

 		\item Step 2: Let $y=\pi_{x}-(g-1)$. 
 		Find all the non-$g$-gap-$\ell$-level excedance letters of $\sigma$  located     weakly to the right of   $\sigma_y$  and  strictly to the left $\sigma_{\ell}$, say $\sigma_{i_1}, \sigma_{i_2}, \ldots, \sigma_{i_k}$ with  $\ell>i_1>i_2>\ldots>i_k\geq y$. Equivalently, we have $$
 		\{i\mid \sigma_{i}<\ell\,\, \mbox{or}\,\, \sigma_{i}<i+g,\,\, y\leq i<\ell\}=\{i_1, i_2, \ldots, i_k\}.
 		$$
 		 Find all the  $g$-gap $\ell$-level excedance letters of $\sigma$ that occur    to  the left of $\sigma_y$ and  are less than   $e_t$  but not less than $\pi_x$, say $\sigma_{j_1}, \sigma_{j_2}, \ldots, \sigma_{j_a}$ with $ j_1<j_2<\cdots<j_a<y$ and $j_p=x$.
 		 Equivalently, we have $$
 		 \{i\mid \sigma_i\in g\Excl_{\ell}(\sigma),\,\, \pi_x\leq \sigma_i<e_t, \,\ i<y\}=\{j_1, j_2, \ldots, j_a\}.
 		 $$
 		 Generate a sequence  $\lambda_1\lambda_2\ldots \lambda_a$ from $\sigma_{j_1}\sigma_{j_2}\ldots \sigma_{j_a}$ by removing $\sigma_{j_p}=\sigma_x$ and inserting $e_t$ as the $(a+1-p)$-th letter (counting from left  to right) of the resulting sequence. 
 	 	Then construct a permutation $\theta$   from $\sigma$ as follows.
 	 	\begin{itemize}
 	 		\item Replace     $\sigma_{i_j}$ with $\sigma_{i_{j+1}}$ for all $1\leq j\leq k$  with the convention that $i_{k+1}=x$ 	 and insert $ \sigma_{i_1}$ immediately before $\sigma_{\ell}$. 
 	 		 	 		\item Replace $\sigma_{j_i}$ with $\lambda_i$ for all $1\leq i\leq a$.
 	 		
 	 	\end{itemize}
  \item  Step 3:	Find all the  $g$-gap $(g+\ell)$-level excedance letters of $\theta$  located  to  the left of $\theta_y$,   say $\theta_{k_1}, \theta_{k_2}, \ldots, \theta_{k_b}$ where  $  k_1<k_2<\cdots<k_b<y$ and $\theta_{k_q}=e_t$ for some $1\leq q\leq b$.  Generate a sequence $\mu_1\mu_2\ldots \mu_b$ from $\theta_{k_1}\theta_{k_2}\ldots \theta_{k_b}$
  by removing $ e_t$ and inserting $e_t$ as the $(b+1-q)$-th letter (counting from left to right) of the resulting sequence.
     Then construct  a permutation $\tau=\tau_1\tau_2\ldots \tau_n$ from $\theta$ by  replacing    $\theta_{k_i}$ with $\mu_i$ for all $1\leq i\leq b$.
   \item Step 4: Find all the $g$-gap $ (g+\ell) $-level excedance letters of $\tau$ that occur   to the left of $\tau_\ell$, say $e'_1, e'_2, \ldots, e'_m$ with  $e'_1<e'_2<\cdots< e'_m$. Define   $ \gamma_{n, g, \ell}(\pi,c)$ to be the permutation obtained from $\tau$  by   replacing 	$e'_j$ with $e'_{j-1}$ for all $1< j\leq m$ and replacing $e'_1$ with $e'_m$.

 	\end{itemize}

 	 \end{framed}
 Take $\pi= 7962(10)48315$, $g=3$ and $\ell=6$ for example.  Clearly,  we have $e_1=6$, $e_2=7$,  $e_3=9$,  $e_4=10$, $3\Excl_9(\pi)=\{9, 10\}$, $\mathcal{A}(\pi)=\{2,5\}$, $\mathcal{B}(\pi)=\{4\}$  and $\mathcal{C}(\pi)=\{1,3\}$.
 Then the $3$-gap-$6$-level-den-labeling    of the index sequence of  $\pi$ is given by$$
 	\begin{array}{llllllllllll} 
 		&_{5}1&_{8}2&_{6}3&_{7}4&_{9}5&_{10}{\bold 6}&_{4}{\bold 7}&_{3}  8&_{2}{\bold 9}&_{1}{\bold 10}&_{0}. 
 		
 	\end{array}
 	$$  It is easily seen that the space labeled by $6$ in the   $3$-gap-$6$-level-den-labeling of the index sequence of $\pi$ is immediately before the index $x=3$ and $x\in \mathcal{C}(\pi)$.  Clearly, we have $y=\pi_3-(g-1)=4$. 
 	 	 First we generate a permutation $\sigma=7(10)62(11)48315$ from $\pi$ by replacing $e_3=9$ with $e_4=10$ and replacing $e_4=10$ with $e_5=11$.  It is easily seen that the non-$3$-gap-$6$-level excedance letter  located weakly to the right of $\sigma_y=\sigma_4=2$ and strictly to the left of $\sigma_{\ell}=\sigma_6=4$ is given by $\sigma_4=2$. Moreover, the $3$-gap $6$-level excedance letters of $\sigma$ that occur to the left of $\sigma_y=\sigma_4=2$ and are less than $e_3=9$ but not less than $\pi_x=\pi_3=6$ are given by $\sigma_{1}=7$ and $\sigma_3=6$.
 	 	   Generate a sequence $\lambda=\lambda_1\lambda_2=97 $  from $\sigma_1\sigma_3=76$ by removing $\sigma_3=6$ and inserting $e_3=9$ as the leftmost letter of the resulting sequence. 	 	 
 	 	  Then construct a permutation $\theta= 9(10)76(11)248315$ from $\sigma$  by the following procedure.
 	 	 \begin{itemize}
 	 	 	\item  Replace $\sigma_4=2$ with $\sigma_3=6$ and insert $\sigma_4=2$ immediately before $\sigma_6=4$.
 	 	  	 	 	\item Replace $\sigma_1$ with $\lambda_1=9$ and replace $\sigma_3=6$ with $\lambda_2=7$. 
 	 	 \end{itemize}
 It is easily checked  that 	the $3$-gap $9$-level excedance letters of $\theta$ that occur to the left of $\theta_y=\theta_4=6$ are given by $\theta_1=9$ and $\theta_2=10$.	 
 Generate a sequence $\mu=\mu_1\mu_2=(10)9$ from $\theta_1\theta_2=9(10)$ by removing $\theta_1=e_3=9$ and inserting it as    the rightmost letter of the resulting sequence.  Then we construct a permutation $\tau=(10)976(11)248315$ from $\theta$ by replacing  $\theta_1=9$ with $\mu_1=10$,  and replacing
 	$\theta_2=10$ with $\mu_2=9$.  
 Clearly, the $3$-gap 	 $9$-level excedance letters of $\tau$ that occur to the left of $\tau_\ell=\tau_6=2$ are given by   $e'_1=\tau_2=9$, $e'_2=\tau_1=10$ and $e'_3=\tau_5=11$.	 Finally, we get a permutation  $\gamma_{11, 3,6}(\pi,6)=9(11)76(10)248315  $ from $\tau$  by replacing $e'_2=10$ with $e'_1=9$,    replacing $e'_3=11$ with $e'_2=10$, and  replacing $e'_1=9$ with $e'_3=11$.

\begin{lemma}\label{lemgamma1}
	Let $g, \ell\geq 1$ and $n\geq g+\ell$.  For any
	$(\pi, c)\in \mathcal{Z}_{n,g, \ell} $,   we have
	\begin{equation}\label{gammare1}
		g\exc_\ell(\gamma_{n,g, \ell}(\pi, c))= g\exc_\ell(\pi)		 	 
	\end{equation}
	and 
	\begin{equation}\label{gammare2}
		g\den_\ell(\gamma_{n,g, \ell}(\pi, c))=g\den_\ell(\pi)+c.
	\end{equation}
\end{lemma}

The following lemma will play an essential role in the proof of Lemma \ref{lemgamma1}. 
\begin{lemma}\label{lemgamma2}
	Fix $g, \ell\geq 1$. Let $\tau=\tau_1\tau_2\ldots \tau_n\in \mathfrak{S}_n $. Assume that $\tau_k=x\geq g$ with  $k\leq x\leq k+(g-1)$. Suppose  that there are exactly $p$  $g$-gap $\ell$-level  excedance letters that occur    to the left of $\tau_{x-(g-1)}$  and  are  greater  than $x$. 
	If $i\in g\Excp_{\ell}(\tau) $ for all $x-(g-1)\leq i<k$,  then  we have
	\begin{equation}\label{mainlem1}
		\begin{array}{lll}
	 &|\{i\mid \tau_i>x,\,\, \tau_i\notin g\Excl_{\ell}(\tau),\,\, i<k\}| &=0
	 \end{array}
	 \end{equation}
 and
 \begin{equation}
 \begin{array}{lll}\label{mainlem2}
	 	  & |\{i\mid \tau_i<x,\,\, \tau_i\notin g\Excl_{\ell}(\tau),\,\, i>k\}|&=p+g-1.
	  \end{array}
\end{equation}

\end{lemma}
\pf  
For any letter $y>x$  with $\tau_z=y$ for some $z<x-(g-1)$, we have $\tau_z=y>x>z+(g-1)$ and hence $y$ is a $g$-gap $\ell$-level   excedance letter. This implies that     there does not exist any non-$g$-gap-$\ell$-level  excedance letters that  occur to the left  of $\tau_{x-(g-1)}$ and   are greater than $x$.  Then (\ref{mainlem1}) follows from the fact that  $i\in g\Excp_\ell(\tau)$   for all $x-(g-1)\leq i<k$.

Assume   that there are exactly $q$   letters that occur to the left of $\tau_k=x$  and are less than $x$ but not less than $k$. This yields that there are   exactly $ x-k-q$   non-$g$-gap-$\ell$-level  excedance letters  that occur  to the right of $\tau_k=x$ and  are less than $x$ but not less than $k$. 
 As  $i\in g\Excp_{\ell}(\tau) $ for all $x-(g-1)\leq i<k$, we have $\tau_i>x$ for all $x-(g-1)\leq i<k$. Then by  (\ref{mainlem1}), there are exactly $p+k-x+(g-1) $  letters    that are greater than $x$ and  occur  to the left of $\tau_k=x$. 
 Hence,   the number of letters that occur  to the left of $\tau_k=x$ and  are  greater than or equal to $k$     is given by $p+q+k-x+(g-1) $.  This implies that there are exactly $p+q+k-x+(g-1) $  letters that occur to the right of $\tau_k=x$  and  are  less than $k$. Clearly, such letters are non-$g$-gap-$\ell$-level  excedance letters.   Recall that there are    exactly $ x-k-q$   non-$g$-gap-$\ell$-level  excedance letters  that occur  to the right of $\tau_k=x$ and are less than $x$ but not less than $k$.
   Hence, the total number of  non-$g$-gap-$\ell$-level  excedance letters that occur to the right of $\tau_k=x$ and  are less than $x$ is given by $p+(g-1)$, completing the proof of  (\ref{mainlem2}).
\qed

Now we are in the position to complete the proof of Lemma \ref{lemgamma1}. 

\noindent{\bf Proof of Lemma \ref{lemgamma1}.}  Here we use the  notations of the definition of  the map $\gamma_{n,g, \ell}$. Let $\gamma_{n,g, \ell}(\pi, c)=\tau'=\tau'_1\tau'_2\ldots \tau'_n$. By the same reasoning as in the proof of Lemma \ref{lemalpha1}, one can easily check   that  $ 
g\den_\ell(\sigma)=g\den_\ell(\pi)
$ and $g\exc_\ell(\sigma)=g\exc_\ell(\pi)$. 
Our next goal is to  verify the following facts.

\noindent{\bf Fact ${\bf1''}$.} 
The set of the non-$g$-gap-$\ell$-level excedance letters of $\theta$ located    weakly to the right of    $\theta_y$  and  weakly to the left of   $\theta_{\ell}$  is given by $\{\theta_{i_j}\mid 0\leq j\leq k \}$,  where $\theta_{i_j}=\sigma_{i_{j+1}}$  for all $0\leq j\leq k$ with the convention that $i_{0}=\ell$ and $i_{k+1}=x$. \\
\noindent{\bf Fact ${\bf2''}$.}  The set of the $g$-gap $\ell$-level excedance letters of $\theta$ that occur  to the left of $\theta_y$ and  are greater than $\pi_x$ but not greater than $e_t$ is given by $\{\theta_{j_i}\mid 1\leq i\leq a \}$.
Equivalently, we have $$
\{i\mid \theta_i\in g\Excl_\ell(\theta),\,\, \pi_x<\theta_i\leq e_t, \,\ i<y\}=\{j_1, j_2, \ldots, j_a\}.$$
\noindent{\bf Fact ${\bf3''}$.} For all $i<\ell$,    the letter $\sigma_i$ is a $g$-gap  $\ell$-level excedance letter of $\sigma$ if and only if $\theta_i$ is a $g$-gap  $\ell$-level excedance letter of $\theta$. Moreover, the letter  $\theta_\ell$ is   a non-$g$-gap  $\ell$-level excedance letter of $\theta$.\\ 
\noindent{\bf Fact ${\bf4''}$.} For all $i\geq \ell$, the letter $\sigma_i$ is a $g$-gap  $\ell$-level excedance letter of $\sigma$ if and only if $\theta_{i+1}=\sigma_i$ is a $g$-gap  $\ell$-level excedance letter  of $\theta$.\\  
\noindent{\bf Fact ${\bf5''}$.} $g\Excl_{g+\ell}(\theta)= \{e_i\mid t\leq i\leq s+1\}$ with the convention that $e_{s+1}=n$.\\
Fact $4''$ could be verified by the same reasoning as in the proof of Fact $4$ and the proof is omitted here. 
Clearly,  the procedure from $\pi$ to $\sigma$ leaves  all the letters not greater than $g+\ell$ in their place.  Then $\pi_x<\ell+g$  implies that $\pi_x=\sigma_x$.
 By the construction of $\theta$, we have $\theta_{i_j}=\sigma_{i_{j+1}}$ for all $1\leq j\leq k$ with the convention that $i_{k+1}=x$ and $i_{0}=\ell$.  Since
$y=\pi_x-(g-1)$ and  $\sigma_x=\pi_x$,  we have $\theta_{i_k}=\sigma_{x}=\pi_x=y+(g-1)\leq i_k+(g-1)$. This yields that $\theta_{i_k}$ is not a $g$-gap $\ell$-level excedance letter of $\theta$. 
Note  that the index $i_j$ is not a $g$-gap $\ell$-level excedance place of $\sigma$ for all $1\leq j\leq k$. Namely, we have  either $\sigma_{i_j}<\ell$ or $  \sigma_{i_j}\leq i_j+(g-1) $ for all $1\leq j\leq k$.   Hence, it follows that
we have either $\theta_{i_{j}}=\sigma_{i_{j+1}}<\ell$ or   $\theta_{i_j}=\sigma_{i_{j+1}}\leq i_{j+1}+(g-1)<i_j+(g-1) $
for all $0\leq j< k$. 
Therefore, we conclude that $\theta_{i_j}$ is a non-$g$-gap-$\ell$-level excedance letter of $\theta$ for all $0\leq j\leq k$.  This completes the proof of Fact $1''$. 

Observe that the sequence $\lambda_1\lambda_2\ldots \lambda_a$ is a permutation of the elements of the set $\{\sigma_{j_i}\mid 1\leq i\leq a,\,\, i\neq p\}\cup \{e_t\}$.
Recall that
$$
\{i\mid \sigma_i\in g\Excl_{\ell}(\sigma),\,\, \pi_x\leq \sigma_i<e_t, \,\ i<y\}=\{j_1, j_2, \ldots, j_a\},
$$
and $\sigma_{j_p}=\sigma_x=\pi_x$.
Then we have $\pi_x<\lambda_i\leq e_t $. 
By the construction of  $\theta$, we have $\theta_{j_i}=\lambda_i$ for all $1\leq i\leq a$.
 Recall that  $y=\pi_x-(g-1)$ and $j_1<j_2<\cdots<j_a<y$. Then, we deduce that 
 $e_t\geq \theta_{j_i}=\lambda_i>\pi_x\geq j_i+g$ for all $1\leq i\leq a$.  This yields that, for all $1\leq i\leq a$, each letter $\theta_{j_i}$ is a   $g$-gap $\ell$-level excedance letter   satisfying that  $\pi_x<\theta_{j_i}\leq e_t$.  This completes the proof of Fact $2''$. 

By the definition  of the procedure from $\sigma$ to $\theta$,   
we have $\theta_z=\sigma_z$ if  $z<\ell$ and $z\notin \{i_j\mid 1\leq j\leq k+1\}\bigcup \{j_i\mid 1\leq i\leq a\}$. Then Fact $3''$ follows directly from Fact $1''$ and Fact $2''$.

 Again by the construction of  $\theta$,  we have $\theta_{i_k}=\sigma_x=\pi_x$. Then, by Fact $1''$ through Fact $4''$,     one can easily check that the procedure from $\sigma$ to $\theta$    transforms  $e_t$ into a $g$-gap $\ell$-level excedance letter,   transforms  $\sigma_x=\pi_x$ into a non-$g$-gap-$\ell$-level excedance letter,   and preserves all the other  $g$-gap $\ell$-level excedance letters of $\sigma$. By the definition of the procedure from $\pi$ to $\sigma$, it is easily checked  that $g\Excl_{g+\ell}( \sigma)=\{ e_i\mid t<i\leq s\} $.  Therefore,  we deduce that $g\Excl_{g+\ell}( \theta)=\{ e_i\mid t\leq i\leq s\}$ as desired, completing the proof of Fact~$5''$.
So far, we have completed the proof   of Fact $1''$ through Fact $5''$.

It is apparent that   $\theta_i=\tau_i=\tau'_i$  for all $i\geq \ell$.  By Fact $3''$  and  Fact $4''$,  it is easily seen that the procedure from $\sigma$ to $\theta$   preserves the number of  $g$-gap $\ell$-level excedance letters that occur weakly to the right of $\sigma_\ell$. Hence, we derive that 
$$
g\exc_\ell(\tau') = g\exc_\ell(\theta)=g\exc_\ell(\sigma)= g\exc_\ell(\pi), 
$$ 
completing the proof of (\ref{gammare1}).

 Now we proceed to prove (\ref{gammare2}).   
Recall that 
$$
g\den_\ell(\sigma)=\sum\limits_{i\in g\Excp_\ell(\sigma)} (i+g-1)+{\inv}(g\Exc_{\ell}(\sigma))+{\inv}(g\Nexc_{\ell}(\sigma)). 
$$
By   Fact $3''$  and Fact $4''$,   it is easily seen that the procedure from $\sigma$ to $\theta$  will increase the sum $\sum\limits_{i\in g\Excp_\ell(\sigma)}(i+g-1)$ by $g\exc_\ell(\sigma)$.  That is, we have
\begin{equation}\label{gammare3}
	\sum\limits_{i\in g\Excp_\ell(\theta)}(i+g-1)=g\exc_\ell(\sigma)+\sum\limits_{i\in g\Excp_\ell(\sigma)}(i+g-1).
\end{equation}

  Assume that 
  $$
  |\{i\mid  \pi_i\in g\Excl_{\ell}(\pi), \,\, \pi_i\geq e_t,   \,\, i<y\}|=u,\,\,\,\,\, |\{i\mid  \pi_i\in g\Excl_{\ell}(\pi), \,\, \pi_i\geq e_t, \,\,  y\leq i<\ell\}|=u',
  $$
    $$
  |\{i\mid \pi_i\in g\Excl_{\ell}(\pi), \,\,\pi_x\leq \pi_i<e_t,  \,\,  i<y\}|=v,
  $$
  and 
  $$
  |\{i\mid \pi_i\in g\Excl_{\ell}(\pi), \,\,\pi_x\leq \pi_i<e_t, \,\,  y\leq i<\ell\}|=v'.
  $$
  Then according to the rules specified in the $g$-gap-$\ell$-level-den-labeling, we have $$c=(g-1)+g\exc_{\ell}(\pi)+u+u'+v+v' $$  since there are exactly $g\exc_{\ell}(\pi)+u+u'+v+v' $ $g$-gap $\ell$-level excedance letters of $\pi$ that are greater than or equal to $\pi_x$. 
  Recall that  the procedure from   $\pi$ to $\sigma$    replaces  a $g$-gap $(g+\ell)$-level excedance letter $e_i$ with $e_{i+1}$ for all $t\leq i\leq s$ with the convention $e_{s+1}=n$ and leaves the other letters fixed. Then, it is easily checked that
  $$
  |\{i\mid \sigma_i\in g\Excl_{\ell}(\sigma), \,\, \sigma_i> e_t,   \,\, i<y\}|=u, \,\,\,\,\,\,   |\{i\mid \sigma_i\in g\Excl_{\ell}(\sigma), \,\, \sigma_i> e_t, \,\,  y\leq i<\ell\}|=u', 
  $$
     $$
  |\{i\mid \sigma_i\in g\Excl_{\ell}(\sigma), \,\,\pi_x\leq \sigma_i< e_t,  \,\,  i<y\}|=v,
  $$
  and 
  $$
  |\{i\mid \sigma_i\in g\Excl_{\ell}(\sigma), \,\,\pi_x\leq \sigma_i<e_t, \,\,  y\leq i<\ell\}|=v'.
  $$
  
  Recall that  we have $\theta_z=\sigma_z$ if  $z<\ell$ and $z\notin \{i_j\mid 1\leq j\leq k\}\bigcup \{j_i\mid 1\leq i\leq a\}$. 
  Then by Fact $1''$ and Fact $2''$, we derive   that $$
  |\{i\mid \theta_i\in g\Excl_{\ell}(\theta), \,\,\theta_i> e_t,   \,\, i<y\}|=u,\,\,\,\,\,\,\, |\{i\mid \theta_i\in g\Excl_{\ell}(\theta), \,\, \theta_i> e_t, \,\,   y\leq i<\ell\}|=u',
  $$
    $$
  |\{i\mid \theta_i\in g\Excl_{\ell}(\theta), \,\,\pi_x< \theta_i\leq e_t,  \,\,  i<y\}|=v,
  $$
  and 
  $$
  |\{i\mid \theta_i\in g\Excl_{\ell}(\theta), \,\,\pi_x<\theta_i\leq e_t, \,\,  y\leq i<\ell\}|=v'.
  $$
  Then there are exactly $u+v$ $g$-gap $\ell$-level excedance letters that are greater than $\pi_x$ and occur to the left of $\theta_y=\theta_{\pi_x-(g-1)}$ in $\theta$. 
  Recall that $\theta_{i_k}=\sigma_x=\pi_x$ with $i_k\leq \pi_x =y+(g-1)\leq   i_k+(g-1)$ and $\pi_x\geq x+g> g$.  Moreover, by Fact $1''$, we have $j\in g\Excp_{\ell}(\theta)$ for all $y\leq j<i_k$.  Then, by Lemma \ref{lemgamma2},   the replacement of $\sigma_{i_k}$ with $\sigma_x=\pi_x$ in the procedure from $\sigma$ to $\theta$ will increase $\inv(g\Nexc_{\ell}(\sigma))$ by  $u+v+g-1$.   That is, we have  \begin{equation}\label{gammare4}
     	\inv(g\Nexc_{\ell}(\theta))=\inv(g\Nexc_{\ell}(\sigma))+u+v+g-1.
     \end{equation}
 
 Recall that $$
\{i\mid \sigma_i\in g\Excl_{\ell}(\sigma),\,\, \pi_x\leq \sigma_i<e_t, \,\ i<y\}=\{j_1, j_2, \ldots, j_a\}, 
$$
and  the sequence $\lambda_1\lambda_2\ldots \lambda_a$ is obtained  from $\sigma_{j_1}\sigma_{j_2}\ldots \sigma_{j_a}$ by removing $\sigma_{j_p}=\sigma_x=\pi_x$ and inserting $e_t$ as the $(a+1-p)$-th letter (counting from  left to right) of the resulting sequence.  Clearly, we have 
$\inv(\lambda_1\lambda_2\ldots \lambda_{a})=\inv(\sigma_{j_1}\sigma_{j_2}\ldots \sigma_{j_a})$. 
 Then  the replacement of    $\sigma_{j_i}$ with $\lambda_i$ for all $1\leq i\leq a$ will increase $\inv(g\Exc_{\ell}(\sigma))$ by  the number of $g$-gap $\ell$-level excedance letters that are  less than $e_t$ but not less than $\pi_x$  and occur weakly to the right of  $\sigma_y$ in $\sigma$.  Recall that $$
 |\{i\mid \sigma_i\in g\Excl_{\ell}(\sigma), \,\,\pi_x\leq \sigma_i<e_t, \,\,  y\leq i<\ell\}|=v'.
 $$  Hence, we conclude that
  \begin{equation}\label{gammare5}
  	\inv(g\Exc_{\ell}(\theta))=\inv(g\Exc_{\ell}(\sigma))+v'.
  \end{equation}
  
 Invoking (\ref{gammare3})-(\ref{gammare5}), we deduce that
\begin{equation}\label{gammare6}
\begin{array}{ll}
	g\den_\ell(\theta)&=\sum\limits_{i\in g\Excp_\ell(\theta)} (i+g-1)+{\inv}(g\Exc_{\ell}(\theta))+{\inv}(g\Nexc_{\ell}(\theta))\\
	&=	g\den_\ell(\sigma) +(g-1)+g\exc_{\ell}(\sigma)+u+v+v'\\ 
	&=  g\den_\ell(\pi) +(g-1)+g\exc_{\ell}(\pi)+u+v+v',
	\end{array}
\end{equation}
where the last equality follows from the fact that $g\exc_\ell(\pi)=g\exc_\ell(\sigma)$ and 
$g\den_\ell(\pi)=g\den_\ell(\sigma)$.

 Recall that  $\mu_1\mu_2\ldots \mu_b$ is obtained from $\theta_{k_1}\theta_{k_2}\ldots \theta_{k_b}$ by removing $\theta_{k_q}=e_t$ and inserting $e_t$ as its $(b+1-q)$-th letter (counting from left  to right) of the resulting sequence.  Since
  each  $\theta_{k_i}$ is a $g$-gap $(g+\ell)$-level excedance letter of $\theta$, we have $\mu_i\geq g+\ell$ for all $1\leq i\leq b$.   Hence, the replacement of $\theta_{k_i}$ with $\mu_i$ in the procedure from $\theta$ to $\tau$   preserves the set of $g$-gap $\ell$-level excedance places  and leaves    the non-$g$-gap-$\ell$-level excedance letters in their place. Hence, this procedure  
  does not affect the sum   $\sum\limits_{i\in g\Excp_\ell(\theta)} (i+g-1)$ and    ${\inv}(g\Nexc_{\ell}(\theta))$.   
  Since $e_t$ is the smallest $g$-gap $(g+\ell)$-level excedance letters of $\theta$ by Fact $5''$, it follows that
    $$\inv(\mu_1\mu_2\ldots \mu_b)=\inv(\theta_{k_1}\theta_{k_2}\ldots \theta_{k_b})+(b-q)-(q-1).$$  Then, replacing  
 $\theta_{k_i}$ with $\mu_i$ for all $1\leq  i\leq b$ in the procedure from $\theta$ to $\tau$ will increase ${\inv}(g\Exc_{\ell}(\theta))$ by $b-q$ and decrease  ${\inv}(g\Exc_{\ell}(\theta))$ by $q-1$. 
 Hence, we derive that
 \begin{equation}\label{gammare8}
 	\begin{array}{ll}
 		g\den_\ell(\tau)&=\sum\limits_{i\in g\Excp_\ell(\tau)} (i+g-1)+{\inv}(g\Exc_{\ell}(\tau))+{\inv}(g\Nexc_{\ell}(\tau))\\
 		&=	g\den_\ell(\theta) -(q-1)+(b-q).\\ 
 	\end{array}
 \end{equation}

 Since  $e'_j\geq g+\ell$ for all $1\leq j\leq m$, the procedure from $\tau$ to $\tau'$   preserves the set of $g$-gap $\ell$-level excedance places and leaves all the non-$g$-gap-$\ell$-level excedance letters in their place. This implies that the  procedure from $\tau$ to $\tau'$ does not affect  the sum  $\sum\limits_{i\in g\Excp_\ell(\tau)} (i+g-1)$ and    ${\inv}(g\Nexc_{\ell}(\tau))$.
Observe that  $e'_m$ is the greatest $g$-gap $(g+\ell)$-level excedance letters located to the left of $\tau$ and   $e'_1=e_t$ is the $(b+1-q)$-th  $g$-gap $(g+\ell)$-level excedance letter of $\tau$ (counting from left to right).  Therefore,      replacing $e'_1=e_t$ with $e'_m$ in the procedure from $\tau$ to $\tau'$  will increase  ${\inv}(g\Exc_{\ell}(\tau))$ by $m-(b+1-q)$ and decrease ${\inv}(g\Exc_{\ell}(\tau))$ by $b-q$. Therefore, we derive that  
 \begin{equation}\label{gammare9}
 		\begin{array}{ll}
 		g\den_\ell(\tau')&=\sum\limits_{i\in g\Excp_\ell(\tau')} (i+g-1)+{\inv}(g\Exc_{\ell}(\tau'))+{\inv}(g\Nexc_{\ell}(\tau'))\\
 		&=	g\den_\ell(\tau) +m-(b+1-q)-(b-q). \\ 
 	\end{array}
 \end{equation}

 Invoking (\ref{gammare6})-(\ref{gammare9}), we derive that 
 \begin{equation}\label{gammare10}
 g\den_\ell(\tau')= g\den_\ell(\pi) +(g-1)+g\exc_{\ell}(\pi)+u+v+v'+(m-b).
 \end{equation} 

 Clearly, the procedure from $\theta$ to $\tau$ preserves the set of  the $g$-gap $(g+\ell)$-level excedance places of $\theta$.  Hence, there are exactly $m$  $g$-gap $(g+\ell)$-level excedance letters of $\theta$ located to the left of $\theta_\ell$.  Note that there are exactly $b$ $g$-gap $(g+\ell)$-level excedance letters of $\theta$ located to the left of $\theta_y$. Then, the number  of the  $g$-gap $(g+\ell)$-level excedance letters of $\theta$ located weakly  to the right of $\theta_y$ and strictly to the left of $\theta_{\ell}$ is given by $m-b$. On the other hand, we have  $$
|\{i\mid \theta_i\in g\Excl_{\ell}(\theta), \,\, \theta_i> e_t, \,\,  y\leq i<\ell\}|=u'.
$$
This implies that $m-b=u'$. Hence,  from  (\ref{gammare10}), it follows  that 
$$
g\den_\ell(\tau')= g\den_\ell(\pi) +(g-1)+g\exc_{\ell}(\pi)+u+v+v'+u'=g\den_\ell(\pi)+c
$$
as desired, completing the proof of  (\ref{gammare2}).

 \qed
 
  \begin{lemma}\label{lemgamma3}
 	Let $g, \ell\geq 1$ and $n\geq g+\ell$.  We have $\gamma_{n,g, \ell}(\pi, c)\in  \mathfrak{S}^{(3)}_{n, g, \ell} $. 
 \end{lemma}
 \pf   Here we use the  notations of the definition of  the map $\gamma_{n,g, \ell}$.  Let $\gamma_{n,g, \ell}(\pi, c)=\tau'=\tau'_1\tau'_2\ldots \tau'_n$. 
  First we aim to show that  $\tau'\in  \mathfrak{S}_{n, g, \ell}-\mathfrak{S}^{(1)}_{n, g, \ell} $. 
 Clearly, we have $\tau'_{\ell}=\tau_{\ell}=\theta_{\ell}$.
 By Fact $3''$, the index $\ell$ is not a $g$-gap $\ell$-level excedance place of $\theta$.
 Hence, we have $\tau'_{\ell}<g+\ell$ and thus $\tau'\notin \mathfrak{S}^{(1)}_{n, g, \ell}$. In order to show that $\tau'\in  \mathfrak{S}_{n, g, \ell}-\mathfrak{S}^{(1)}_{n, g, \ell} $, it remains to  the index $z$ is not a $g$-gap $(g+\ell)$-level  grande-fixed place for any $z\geq \ell$.   
 Recall that $\tau'_i=\tau_{i}=\theta_i$ for all $i\geq \ell$. In order to prove $\tau'\in  \mathfrak{S}_{n, g, \ell}-\mathfrak{S}^{(1)}_{n, g, \ell} $,  it suffices to show that 
 the index $z$ is not a $g$-gap   $(g+\ell)$-level grande-fixed place of $\theta$ for any $z\geq  \ell$.   By the same reasoning as in the proof of    Lemma~\ref{lemalpha2}, one can verify that  the index $z$ is not a $g$-gap $(g+\ell)$-level   grande-fixed place of $\theta$ for any $z>\ell$ by employing     Fact $5''$,   and the proof is omitted here.  Note that $\theta_{\ell}$ is a non-$g$-gap-$\ell$-level excedance letter of $\theta$ by Fact $3''$. This implies that $\theta_{\ell}<g+\ell$ and thus
 the index $\ell$ is  not a $g$-gap $(g+\ell)$-level   grande-fixed place of $\theta$.  Hence, we conclude that the index $z$ is not a $g$-gap   $(g+\ell)$-level grande-fixed place of $\theta$ for any $z\geq  \ell$, and thus $\tau'\in \mathfrak{S}_{n, g, \ell}-\mathfrak{S}^{(1)}_{n, g, \ell} $.
 
  In order to prove that $\tau'\in  \mathfrak{S}^{(3)}_{n, g, \ell} $,  it remains to show that $\tau'$ is improper.  
  To this end,  we shall prove that  the resulting permutation $\tau'$  verifies the following desired properties.      
 \begin{itemize}
 
 	\item[\upshape(\rmnum{1})]  $e'_m$ is the greatest  $g$-gap $(g+\ell)$-level excedance letter located  to the left of $\tau'_\ell$ in $\tau'$.		
 	\item[\upshape(\rmnum{2})]  Assume that $\tau'_{z}=e'_m$ and the set of 
 	non-$g$-gap-$\ell$-level excedance letters located  between $\tau'_z$  and    $\tau'_\ell$ (including $\tau'_{\ell}$)  is given by 
 	$\{\tau'_{r_j}\mid  0\leq j\leq k'\}$ with  $\ell=r_{0}>r_1>\cdots> r_{k'}>r_{k'+1}=z$. Then we have  $r_{k+1}<y$ and    $\tau'_{r_j}=\sigma_{i_{j+1}}$ for all $0\leq j\leq k$ with the convention that $i_0=\ell$ and $i_{k+1}=x$. 
 \end{itemize}

    Recall that $e'_m$ is the greatest $g$-gap $(g+\ell)$-level excedance letter  located to the left of  $\tau_{\ell}$ in $\tau$. Clearly, the procedure from $\tau$ to $\tau'$ preserves the set of $g$-gap $(g+\ell)$-level excedance letters  located to the left of  $\tau_{\ell}$ in $\tau$. This implies that $e'_m$ is also the greatest $g$-gap $(g+\ell)$-level excedance letter  located to the left of  $\tau'_{\ell}$ in $\tau'$, completing the proof of (\rmnum{1}). 
     
    By the construction of $\tau$, we have $\tau_{k_{b+1-q}}=e'_1=e_t$ and $k_{b+1-q}<y$.  Recall that $\tau'$ is obtained from $\tau$ by replacing $e'_i$ with $e'_{i-1}$ for all $2<i\leq m$ and replacing $e'_1=e_t$ with  $e'_m$.  This yields that  $\tau'_{k_{b+1-q}}=e'_m$ and hence $z=k_{b+1-q}<y$. Observe that both the procedure 
    from $\theta$ to $\tau$ and the procedure from $\tau$ to $\tau'$   preserve the  set of  $g$-gap $\ell$-level excedance places and leave the  non-$g$-gap-$\ell$-level excedance letters in their place.   Then  (\rmnum{2}) follows directly from    Fact $1''$.

  By properties (\rmnum{1}) and  (\rmnum{2}), in order to prove 
  $\tau'$ is improper, it suffices to show that  $\tau'_{i_k}>\ell$ and $\tau'_{r_k}>r_{k+1}+(g-1)$. 
  By  (\rmnum{2}), we have $r_{k+1}<y=\pi_x-(g-1)$ and $\tau'_{r_k}=\sigma_x$.  Recall  that $\sigma_x=\pi_x>\ell$ by the construction of $\sigma$.  Hence, we deduce that  $\tau'_{r_k}>\ell$ and $\tau'_{r_k}=\pi_x>r_{k+1}+(g-1)$ as desired, completing the proof. \qed
     	 
     	 \begin{lemma}\label{lemgamma4}
     	 	Let $g, \ell\geq 1$ and $n\geq g+\ell$.  The map  $\gamma_{n, g,\ell}  $ is injective.
     	 \end{lemma} 
     	\pf   Here we use the  notations of the definition of  the map $\gamma_{n,g, \ell}$.  Let $\gamma_{n,g, \ell}(\pi, c)=\tau'=\tau'_1\tau'_2\ldots \tau'_n$.      	   	 
     	 By (\ref{gammare2}), we have $c=g\den_\ell(\tau')-g\den_\ell(\pi)$. In order to  show that $\gamma_{n,g,\ell}$ is injective, it is sufficient  to show that we can recover $\pi$ from $\tau'$.   To this end, we give a  description of the procedure  to recover the permutation $\pi$ from $\tau'$.
     	  
     	 \begin{itemize}
     	 	\item Step $1$:   Find all the  $g$-gap  $(g+\ell)$-level excedance  letters of $\tau'$  located  to the left of $\tau'_\ell$, say $\bar{e}_1, \bar{e}_2, \ldots, \bar{e}_{m'} $ with $\bar{e}_1<\bar{e}_2<\cdots<\bar{e}_{m'}=\tau'_{z}$. 	Construct a permutation $\tau$ from $\tau'$ by replacing $\bar{e}_{i-1}$ with $\bar{e}_{i}$ for all $1<i\leq m'$ and replacing $\bar{e}_{m'}$ with $\bar{e}_1$. 
     	 		\item Step $2$: 
     	 	Find all the non-$g$-gap-$\ell$-level excedance letters,   say   $ \tau_{c_0},  \tau_{c_1}, \ldots,  \tau_{c_{u}}$ with $z=c_{u+1}<c_{u}< c_{u-1}<\cdots<c_{0}=\ell$,  that are  located  between $\tau_z=\bar{e}_1$  and  $ \tau_\ell$ (including $ \tau_\ell$) in $ \tau$. 
     	 	 Choose $k'$ to be the least integer such that $ \tau_{c_{k'}}\geq c_{k'+1}+g$ and $\tau_{c_{k'}}>\ell$.
     	 	Let $y'=\tau_{c_{k'}}-(g-1)$.  Find all the $g$-gap  $(g+\ell)$-level excedance  letters of $\tau$  located  to the left of $\tau_{y'}$,  say $\tau_{d_1}, \tau_{d_2}, \ldots, \tau_{d_{v}}$ where $d_1<d_2<\cdots<d_{v}<y'$ and  $d_{q'}=z$. Generate a sequence $\mu'_1\mu'_2\ldots \mu'_{v}$ from $\tau_{d_1}\tau_{d_2}\ldots \tau_{d_{v}}$ by removing $\bar{e}_1$ and inserting $\bar{e}_1$ as the $(v+1-q')$-th letter (counting from left to right) of the resulting sequence.  Then construct a permutation $\theta=\theta_1\theta_2\ldots \theta_n$ from $\tau$  by replacing  $\tau_{d_i}$ with $\mu'_i$  for all $1\leq i\leq v$. 
     	 	 
     	 		\item Step $3$:  Find all the $g$-gap $\ell$-level excedance letters that occur  to the left of $\theta_{y'}$ and are greater than $\theta_{c_{k'}}$ but not greater than $\bar{e}_1$, say $\theta_{f_1}$, $\theta_{f_{2}}, \ldots, \theta_{f_{w}}$
     	 	with $f_1<f_2<\cdots < f_{w}<y'$ and $\theta_{f_{p'}}=\bar{e}_1$. Generate a sequence $\lambda'_1\lambda'_2\ldots \lambda'_{w}$ from $\theta_{f_1}\theta_{f_2}\ldots \theta_{f_{w}}$ by removing $\bar{e}_1$ and inserting $\bar{e}_1$ as the $(w+1-p')$-th letter (counting from left to right) of the resulting sequence. Then we construct a permutation $\sigma=\sigma_1  \sigma_2\ldots \sigma_n$ from $\theta$ as follows.
     	 	\begin{itemize}     	 		
     	 	     	 		\item Replace $\theta_{f_i}$ with $\lambda'_i$  for all $1\leq i\leq w$.
     	 		 \item Remove $\theta_{\ell}$, replace $\theta_{c_{j}}$ with $\theta_{c_{j-1}}$
     	 		for all $1\leq j\leq k'$, and replace $\bar{e}_1$ with $\theta_{c_{k'}}$.
     	 		\end{itemize}   	 	
     	 		\item Step $4$:  Find all the  $g$-gap  $(g+\ell)$-level excedance  letters of $\sigma$, say $\bar{e}'_1, \bar{e}'_2, \ldots, \bar{e}'_h$ with $\bar{e}'_1<\bar{e}'_2<\cdots<\bar{e}'_h$.
     	 		Let $\pi$ be the permutation obtained from $\sigma$ by replacing  $\bar{e}'_i$ with $\bar{e}'_{i-1}$ for each $1\leq i\leq h$ with the convention that $\bar{e}'_0=\bar{e}_1$.

     	 \end{itemize} 
       
     	In the following, we aim to show that    after applying  the above procedure to $\tau$  we will recover the permutation $\pi$.  
     	
     	 Since the procedure from $\tau$ to $\tau'$ preserves the set of $g$-gap $(g+\ell)$-level excedance letters located to the left of $\tau_{\ell}$, one can easily check that  the procedure from $\tau'$ to $\tau$  indeed  reverses the  procedure from $\tau$ to $\tau'$.

     	In order to show the procedure from $\tau$ to $\theta$ reverses the procedure from $\theta$ to $\tau$,  it suffices to show that (\rmnum{1}) $\bar{e}_1=e_t$, (\rmnum{2}) $y=y'$, and (\rmnum{3}) $b=v$, $q'=b+1-q$,  and $d_i=k_i$ for all $1\leq i\leq b$.

      By Fact $5''$, it follows that $e_t$ is the smallest $g$-gap $(g+\ell)$-level excedance letter of $\theta$. Recall that $\theta_{j_{a+1-p}}=e_t$ by the construction of $\theta$. Since $j_{a+1-p}<y<\ell$, it follows   $e_t$ is the smallest  $g$-gap $(g+\ell)$-level excedance letter located to  the left of $\theta_{\ell}$  in  $\theta$. Note that $\bar{e}_1$ is the smallest  $g$-gap $(g+\ell)$-level excedance letter   located to the left of $\tau_{\ell}$.  Moreover,   the  procedure from $\theta$ to $\tau$ 
      preserves the set of $g$-gap $(g+\ell)$-level excedance letters   located to the left of $\theta_{\ell}$.  This yields that $e_t=\bar{e}_1$ as claimed in (\rmnum{1}). 
     
     	Note  that $\theta_i=\tau_i$ for all $i\geq y$ by the construction of $\tau$. Then,  by Fact $1''$,  it follows that the set of the non-$g$-gap-$\ell$-level excedance letters located weakly to the right of  $\tau_y$ and  weakly to the left of   $\tau_{\ell}$ is given by  $\{\tau_{i_j}\mid 0\leq j\leq k\}$ where $\ell=i_0>i_1>\cdots> i_k$ and $\tau_{i_j}=\theta_{i_j}=\sigma_{i_{j+1}}$ for all $0\leq j\leq k$ with the convention that $i_{k+1}=x$.  This implies that   $c_{k+1}<y$ and $c_j=i_j$ for all $1\leq j\leq k$. 
     	  
     	  Note that the procedure from $\pi$ to $\sigma$ leaves all the letters less than $g+\ell$ in their place. Then $\pi_x<g+\ell$ would imply that  $\theta_x=\pi_x$. 
     	  Recall that $\sigma_{i_{j}}$ is not a $g$-gap $\ell$-level excedance letter of $\sigma$ for all $1\leq j\leq k$. This implies that  we have  either $\sigma_{i_{j}}<\ell$ or $\sigma_{i_{j}}\leq i_{j}+g-1$ for all $1\leq j\leq k$. 	Hence, we have $\tau_{c_j}=\tau_{i_{j}}=\sigma_{i_{j+1}}<\ell$ or  $\tau_{c_j}=\tau_{i_{j}}=\sigma_{i_{j+1}}\leq i_{j+1}+g-1$ for all $0\leq j<k$. Recall that $\tau_{c_k}=\tau_{i_k}=\sigma_x=\pi_x=y+(g-1)> c_{k+1}+(g-1)$
     	  and $\tau_{c_k}=\tau_{i_k}=\pi_x>\ell$. This  implies that $k$ is the least integer such that $\tau_{c_k}>\ell$ and $\tau_{c_k}\geq c_{k+1}+g$.    Therefore, we  deduce that $k=k'$,  $c_{k'}=c_{k}=i_{k}$ and thus $y=\tau_{i_k}-(g-1)=\tau_{c_{k'}}-(g-1)=y'$ as claimed in (\rmnum{2}).

     	Note that the procedure from $\theta$ to $\tau$ preserves the set of $g$-gap $(g+\ell)$-level excedance letters located to the left of $\theta_y$ as well as the set of $g$-gap $(g+\ell)$-level excedance places located to the left of $\theta_y$. This implies that  $b=v$ and $d_i=k_i$ for all $1\leq i\leq b$. From the definition of the procedure from $\theta$ to $\tau$, we can see that  $\tau_{k_{b+1-q}}=e_t$. This yields that $q'=b+1-q$, completing the proof of (\rmnum{3}). 
     	
     	 Recall that $e_t=\bar{e}_1$,   $k=k'$, $y=y'$ and $i_j=c_j$ for all $0\leq j\leq k$. 
     	In order to show that the procedure from  $\theta$ to $\sigma$ reverses the procedure from $\sigma$ to $\theta$,   it remains to   show that $a=w$, $p'=a+1-p$,   $j_i=f_i$ for all $1\leq i\leq a$.  Fact $2''$ states  that the set of $g$-gap $\ell$-level excedance letters of $\theta$ that are located to the left of $\theta_y$ and are greater than $\pi_x=\theta_{i_k}$ but not greater than $e_t$ is given by $\{\theta_{j_i}\mid 1\leq i\leq a\}$. This implies that $a=w$ and $f_i=j_i$ for all $1\leq i\leq a$. 
     	Recall that $\theta_{j_{a+1-p}}=e_t$. Hence, we have $p'=a+1-p$ as desired.

     	By the construction of $\sigma$, one can easily check that  the set of the $g$-gap $(g+\ell)$-level excedance letters of $\sigma$ is  given by  $\{e_i\mid t<i\leq s+1\}$ with the convention that $e_{s+1}=n$. This ensures that the procedure from $\sigma$ to $\pi$ indeed reverses the procedure  from $\pi$ to $\sigma$   as desired. This completes the proof.  \qed

  Let $\mathcal{W}_{n,g, \ell}$ denote  the subset of pairs $(\pi, c)\in \mathfrak{S}_{n-1}\times \{ 0, 1, \ldots,  n-1\}$ such that either  $c=0$ or  the space labeled by $c$  under the $g$-gap-$\ell$-level-den-labeling of the index sequence of $\pi$    is immediately  before the index $x$   for some $ x\geq \ell$. 
   
   \begin{framed}
   	\begin{center}
   		{\bf The map $\delta_{n,g, \ell}:\mathcal{W}_{n,g, \ell}\longrightarrow \mathfrak{S}_{n}-\mathfrak{S}_{n, g, \ell} $ }
   	\end{center}
   	Let $(\pi, c)\in \mathcal{W}_{n,g, \ell}$.  Suppose   that
       	$$  g\Excl_{\ell}(\pi) =\{e_1, e_{2}, \ldots, e_s\}$$ 
       	and 
       		$$  g\Excl_{g+\ell}(\pi) =\{e_t, e_{t+1}, \ldots, e_s\}$$ 
       	with $e_1<e_2<\cdots<e_s$.  Assume that  when $c>0$,  the space labeled by $c$  under the $g$-gap-$\ell$-level-den-labeling of the index sequence of $\pi$    is  immediately before the index $x$ for some  $x\geq \ell$,  and assume that the $g$-gap $(g+\ell)$-level insertion letter for this space is given by $e_{y}$  for some  $t\leq y\leq s+1$ with the convention that $e_{s+1}=n$. 
   	Define  $\delta_{n,g,\ell}(\pi,c)$ to be the permutation  constructed as follows.
   	\begin{itemize}
   		\item If $c=0$, then let $\delta_{n,g,\ell}(\pi,c)$ be the permutation obtained from $\pi$ by inserting $n$ immediately  after $\pi_{n-1}$.

   		\item For $c\geq 1$, we generate a permutation  $\sigma$ from $\pi$ by replacing $e_j$ with  $e_{j+1}$ for all $y\leq j\leq s$. Then define $\delta_{n,g,\ell}(\pi,c)$ to be the permutation obtained from $\sigma$ by inserting the letter $e_y$ immediately before the letter $\sigma_{x}$.
   		
   	\end{itemize}

   \end{framed}

Continuing with our running example where $\pi= 596284317$, $g=3$ and $\ell=6$,   we have $e_1=6$, $e_2=8$,   $e_3=9$, and $3\Excl_9(\pi)=\{e_3\}$.
Recall that the $3$-gap-$6$-level-den-labeling  and the $3$-gap $9$-level insertion letters   of the index sequence of  $\pi$ are given by$$
\begin{array}{lllllllllll} 
	&_{6}1&_{8}2&_{5}3&_{7}4&_{4}5&_{9}{\bold 6}&_{3}7&_{2}{\bold 8}&_{1}{\bold 9}&_{0}\\
	&  &  &  &  & &\uparrow &\uparrow&\uparrow&\uparrow&\uparrow\\
	&  &  &  &  & &9&9&10&10&10.
	
\end{array}
$$
   It is easily seen that the space labeled by $9$ in the   $3$-gap-$6$-level-den-labeling of the index sequence of $\pi$ is immediately before the index $x=6$ and  the insertion letter for this space is given by $e_3=9$.  First we generate a permutation $\sigma=5(10)6284317$ from $\pi$ by replacing $e_3=9$ with $e_4=10$.   Then we construct   the permutation  $\delta_{10, 3,6}(\pi,9)=5(10)62894317    $   obtained from $\sigma$ by  inserting $e_3=9$ immediately before $\sigma_6=4$.

   Now we proceed to   show that the map $\delta_{n,g, \ell}$ verifies the following desired properties. 
   \begin{lemma}\label{lemdelta1}
   	Let $g, \ell\geq 1$ and $n\geq g+\ell$.  For any
   	$(\pi, c)\in \mathcal{W}_{n,g, \ell} $, we have  
   	\begin{equation}\label{deltare1}
   		g\exc_\ell(\delta_{n,g, \ell}(\pi, c)) =\left\{ \begin{array}{ll}
   			g\exc_{\ell}(\pi)&\, \mathrm{if}\,\, 0\leq c\leq g\exc_{\ell}(\pi)+\ell+g-2,\\
   			g\exc_{\ell}(\pi)+1&\, \mathrm{otherwise}, 
   		\end{array}
   		\right.
   	\end{equation}
   and 
      	\begin{equation}\label{deltare2}
   		g\den_{g+\ell}(\delta_{n,g, \ell}(\pi, c))-g\den_{g+\ell}(\pi)=	g\den_\ell(\delta_{n,g, \ell}(\pi, c))-g\den_\ell(\pi)=c.
   	\end{equation}
    
   \end{lemma}
   \pf Here we use the  notations of the definition of  the map $\delta_{n,g, \ell}$. Let $\delta_{n,g, \ell}(\pi, c)=\tau=\tau_1\tau_2\ldots \tau_n$. It is apparent that the statement holds for $0\leq c\leq g-1$. Now we assume that $c>g-1$. By the same reasoning as in the proof of Lemma \ref{lemalpha1}, one can verify that      
   $g\exc_{\ell}(\pi)=g\exc_{\ell}(\sigma)$,  $g\den_\ell(\sigma)=g\den_\ell(\pi) $   and  $ 
   g\den_{g+\ell}(\sigma)=g\den_{g+\ell}(\pi).$

   In the following, we shall verify the following facts.\\
   \noindent{\bf Fact A:} The letter $\sigma_i$ is a $g$-gap  $\ell$-level excedance letter (resp. $g$-gap  $(g+\ell)$-level excedance letter) of $\sigma$ if and only if $\tau_i$ is a $g$-gap  $\ell$-level excedance letter   (resp. $g$-gap  $(g+\ell)$-level excedance letter) of $\tau$ for all $i<x$. \\ 
    \noindent{\bf Fact B:} For all $i\geq x$, the letter $\sigma_i$ is a   $g$-gap  $(g+\ell)$-level excedance letter of $\sigma$ if and only if $\tau_{i+1}=\sigma_i$ is a   $g$-gap  $(g+\ell)$-level excedance letter   of $\tau$.\\
     \noindent{\bf Fact C:} The letter $\tau_{x}$ is a    non-$g$-gap-$(g+\ell)$-level excedance letter  of $\tau$ if and only if $x=e_y-(g-1)$. \\
   Fact  A  follows directly from the fact that $\tau_i=\sigma_i$ for all $i<x$ and Fact  B  can be verified by the same reasoning as in the proof of Fact 4.  Fact C follows immediately from the fact that  $\tau_x=e_y$.

  Relying on  Fact A through Fact C, one can easily check that  the procedure from $\sigma$ to $\tau$   increases $g\exc_\ell(\sigma)$ by $1$  when  $x<e_y-(g-1)$ and preserves  the number $g\exc_\ell(\sigma)$ when  $x=e_y-(g-1)$. According to the rules specified in the $g$-gap-$\ell$-level-den-labeling the spaces of the index sequence of $\pi$, we have $c=(g-1)+s+1-y   $ and $c\leq  g\exc_{\ell}(\pi)+g+\ell-2$ when $x=e_y-(g-1)$,  and we have $c=(g-1)+s+1-y+x$ and $c>g\exc_{\ell}(\pi)+g+\ell-2$ when $x<e_y-(g-1)$. 
    As $g\exc_\ell(\pi)=g\exc_\ell(\sigma)$, we have
   $ 
   g\exc_\ell(\tau)=g\exc_\ell(\pi)
   $ 
   when $c\leq g\exc_{\ell}(\pi)+g+\ell-2$, and
   $ 
   g\exc_\ell(\tau)=g\exc_\ell(\pi)+1
   $ 
   when $c> g\exc_{\ell}(\pi)+g+\ell-2$, 
   completing the  proof of~\eqref{deltare1}.
   
   Next we proceed to prove ~\eqref{deltare2}.  Assume that there are $p$  $g$-gap $(g+\ell)$-level excedance letters located weakly to the right of $\pi_{x}$ and 
   $q$  $g$-gap $(g+\ell)$-level excedance letters that   occur  to the left of $\pi_{x}$ and  are  not less than $e_y$. 
   According to the rules specified in the  $g$-gap $(g+\ell)$-level insertion letter for the space before the index $x$, we have $e_{y-1}-(g-1)<x\leq e_y-(g-1)$ when $y>t$. This implies that all the $g$-gap $(g+\ell)$-level excedance letters located weakly  to the right of $\pi_{x}$   are  not less than $e_y$. Hence, we have $p+q=s+1-y$. By the construction of $\sigma$, one can easily check that the resulting permutation $\sigma$ verifies the following properties.
   \begin{itemize}
   	\item[\upshape(\rmnum{1})] There are exactly  $p$  $g$-gap $(g+\ell)$-level excedance letters  located weakly to the right of $\sigma_{x}$ and 
   	$q$  $g$-gap $(g+\ell)$-level excedance letters that are located   to the left of $\sigma_{x}$ and are greater than  $e_y$. 
   		\item[\upshape(\rmnum{2})] There does not exist any $g$-gap $(g+\ell)$-level excedance letters  that occur weakly  to the right of $\sigma_{x}$  and  are    less than $e_y$.
   \end{itemize}

   Now we are in the position to prove  ~\eqref{deltare2}    by distinguishing two cases.

        \noindent{\bf Case 1:} $x<e_y-(g-1)$. \\
        Recall that 
   $$
   g\den_\ell(\sigma)=\sum\limits_{i\in g\Excp_\ell(\sigma)} (i+g-1)+{\inv}(g\Exc_{\ell}(\sigma))+{\inv}(g\Nexc_{\ell}(\sigma)), 
   $$   
   and 
   $$
   g\den_{g+\ell}(\sigma)=\sum\limits_{i\in g\Excp_{g+\ell}(\sigma)} (i+g-1)+{\inv}(g\Exc_{g+\ell}(\sigma))+{\inv}(g\Nexc_{g+\ell}(\sigma)). 
   $$   
   
  By  Fact A through  Fact C ,   one can easily check that the procedure from $\sigma$ to $\tau$  does not affect $\inv(g\Nexc_\ell(\sigma))$ and $\inv(g\Nexc_{g+\ell}(\sigma))$. Moreover, Combining Fact A through Fact  C and  properties (\rmnum{1})-(\rmnum{2}),  one can see that this procedure increases the sum $\sum\limits_{i\in g\Excp_\ell(\sigma)}(i+g-1)$ (resp. $\sum\limits_{i\in g\Excp_{g+\ell}(\sigma)}(i+g-1)$) by 
   $x+g-1+p$ and increases $\inv(g\Exc_\ell(\sigma))$ (resp.  $\inv(g\Exc_{g+\ell}(\sigma))$) by $q$.
   Recall that $g\den_\ell(\sigma)=g\den_\ell(\pi)$, $g\den_{g+\ell}(\sigma)=g\den_{g+\ell}(\pi)$, $c=(g-1)+s+1-y+x$ and $p+q=s+1-y$.
      Therefore, we derive  that
      $$
              	g\den_\ell(\tau)-g\den_\ell(\pi)=g\den_\ell(\tau)-g\den_\ell(\sigma)=x+g-1+p+q=c
      $$  
    and 
   $$
   g\den_{g+\ell}(\tau)-g\den_{g+\ell}(\pi)=g\den_{g+\ell}(\tau)-g\den_{g+\ell}(\sigma)=x+g-1+p+q=c
   $$ 
   as claimed in (\ref{deltare2}).

   \noindent{\bf Case 2:} $x=e_y-(g-1)$. \\
     Combining  Fact A through  Fact C and property (\rmnum{1}),    one can easily check that the procedure from $\sigma$ to $\tau$  does not affect $\inv(g\Exc_\ell(\sigma))$ (resp. $\inv(g\Exc_{g+\ell}(\sigma))$) and   increases the sum $\sum\limits_{i\in g\Excp_\ell(\sigma)}(i+g-1)$ (resp. the sum $\sum\limits_{i\in g \Excp_{g+\ell}(\sigma)}(i+g-1)$) by 
   $p$.  By Lemma (\ref{lemgamma2}) and property (\rmnum{1}),  it is apparent that the insertion of $e_y$ immediately before $\sigma_x$ in the procedure from $\sigma$ to $\tau$ increases $\inv(g\Nexc_\ell(\sigma))$ (resp. $\inv(g\Nexc_{g+\ell}(\sigma))$) by $q+(g-1)$. 
   Recall that  $g\den_\ell(\sigma)=g\den_\ell(\pi)$,  $g\den_{g+\ell}(\sigma)=g\den_{g+\ell}(\pi)$,   $c=(g-1)+s+1-y$ and and $p+q=s+1-y$.   
   Therefore, we derive  that
   $$
  g\den_\ell(\tau)-g\den_\ell(\pi)=g\den_\ell(\tau)-g\den_\ell(\sigma)=g-1+p+q=c
  $$  
   and
  $$
  g\den_{g+\ell}(\tau)-g\den_{g+\ell}(\pi)=g\den_{g+\ell}(\tau)-g\den_{g+\ell}(\sigma)=g-1+p+q=c
  $$ 
  as claimed in (\ref{deltare2}).
    This completes the proof. 
   \qed
   
   \begin{lemma}\label{lemdelta2}
  	Let $g, \ell\geq 1$ and $n\geq g+\ell$.  The map  $\delta_{n, g,\ell}  $ is an injection from $\mathcal{W}_{n,g,\ell}$ to $\mathfrak{S}_n-\mathfrak{S}_{n, g, \ell}$. 
  \end{lemma}
  \pf  Here we use the  notations of the definition of  the map $\delta_{n,g, \ell}$. 
  Let $\delta_{n,g, \ell}(\pi, c)=\tau=\tau_1\tau_2\ldots \tau_n$. 
  
\noindent{\bf Claim 1:}	  $x$ is the greatest integer such that  the index $x$ is a   $g$-gap $(g+\ell)$-level   grande-fixed place of $\tau$. \\
    In the following,  we proceed the proof of  Claim 1 by considering two cases.  
  
  \noindent{\bf Case 1:} $c\leq g-1$.\\
  In this case, we have $\tau_x=n$. Clearly, the index $x$ is the greatest integer such that the index $x$ is a $g$-gap $(g+\ell)$-level  grande-fixed place of $\tau$ as desired  in Claim 1.

   \noindent{\bf Case 2:} $c>g-1$.\\     
     In this case, we have $\tau_x=e_y\geq g+\ell$. By the construction of $\sigma$, we have $g\Excl_{g+\ell}(\sigma)=\{e_i\mid t\leq  i\leq s+1\}-\{e_y\}$. Then,  from Fact B and Fact C, it follows that 
  $g\Excl_{g+\ell}(\tau)=\{e_i\mid t\leq i\leq s+1\}$ when $x<e_y-(g-1)$ and $g\Excl_{g+\ell}(\tau)=\{e_i\mid t\leq  i\leq s+1\}-\{e_y\}$ when $x=e_y-(g-1)$.      
  According to the rules specified in the  $g$-gap $(g+\ell)$-level insertion letter for the space before the index $x$,  we have $e_{y-1}-(g-1)<x\leq e_y-(g-1)$ when $y>t$. 
  This yields that if $x<e_y-(g-1)$, then 
   we have  $\tau_x=e_y\geq x+g\geq g+\ell$ and  $[x+g-1, \tau_x)\cap g\Excl_{g+\ell}(\tau)=\emptyset$.  
  Therefore,    we have    either (1) $\tau_x=x+(g-1)$ or (2) $\tau_x\in g\Excl_{g+\ell}(\tau)$ and  $[x+g-1, \tau_x)\cap g\Excl_{g+\ell}(\tau)=\emptyset$, implying that the index $x$ is a $g$-gap $(g+\ell)$-level grande-fixed place of $\tau$. 
  
  In order to prove Claim 1, it  remains to show that there does not  exist any $g$-gap $(g+\ell)$-level  grande-fixed place  $z$ of $\tau$  with  $z>x$.   It is apparent that  $e_{s+1}=n$ is a 
   $g$-gap  $(g+\ell)$-level excedance letter of $\tau$.  This implies that $\tau_i\neq n$ for all $n-(g-1)\leq i\leq n$. Hence, the index $z$ is not a $g$-gap $(g+\ell)$-level  grande-fixed place of $\tau$ for any $z>n-g$. 
   
   Take $z>x$ arbitrarily (if any) such that $\tau_z= z+(g-1)$.  By Fact B, $\tau_z= z+(g-1)$ would imply  that $\sigma_{z-1}\leq z-1+(g-1)$. On the other hand, by the definition of the procedure from $\sigma$ to $\tau$,  we have  $ \sigma_{z-1}=\tau_z=z+(g-1)$,  a contradiction.
   Therefore, we conclude that $\tau_z\neq  z+(g-1)$ for any $z>x$. 
    
  Take $z>x$ arbitrarily (if any) such that $\tau_z\geq  z+g$.  Again by Fact B,  one can see that  $\sigma_{z-1}=\tau_z$ is a $g$-gap  $(g+\ell)$-level excedance letter in $\sigma$. 
  Note that the procedure from $\pi$ to $\sigma$ preserves the set of $g$-gap $(g+\ell)$-level excedance places of $\pi$.  Then the letter $\pi_{z-1}$ is also a $g$-gap  $(g+\ell)$-level excedance letter in $\pi$. Hence we have $\pi_{z-1}\geq z-1+g$  and $\pi_{z-1}=e_{y'} $ for some $t\leq y'\leq s$. 
  Now we proceed to show that $[z+g-1, \tau_z) \cap g\Excl_{g+\ell}(\tau)\neq \emptyset$. We have two cases.
  
  \noindent{\bf Case $1'$:} $x=e_y-(g-1)$.\\
  In this case, we have  
  $e_{y'}=\pi_{z-1}\geq z-1+g\geq x+g=e_y+1\geq e_{y+1}$. This implies that $y'\geq y+1$.
     By the construction of $\sigma$, we have $\sigma_{z-1}>\pi_{z-1}=e_{y'}$.    
   Note that  $g\Excl_{g+\ell}(\tau)=\{e_i\mid t\leq i\leq s+1, i\neq y\}$ when $x=e_y-(g-1)$. Hence, we deduce that 
  $$ 
 e_{y'}\in  [z-1+g, \sigma_{z-1})\cap g\Excl_{g+\ell}(\tau)=[z-1+g, \tau_z)\cap g\Excl_{g+\ell}(\tau)\neq \emptyset. 
  $$

  \noindent{\bf Case $2'$:} $x<e_y-(g-1)$.\\
    Recall that  $x> e_{y-1}-(g-1)$ by the rules specified in the $g$-gap $(g+\ell)$-level insertion letter for the space before the index $x$. Then
  $e_{y'}=\pi_{z-1}\geq z-1+g\geq x+g>e_{y-1}+1\geq e_{y}$. This implies that $y'\geq y$.
  By the construction of $\sigma$, we have $\sigma_{z-1}>\pi_{z-1}=e_{y'}$.    
  Note that   $g\Excl_{g+\ell}(\tau)=\{e_i\mid t\leq i\leq s+1\}$ when $x<e_y-(g-1)$. Hence, we deduce that 
  $$ 
  e_{y'}\in  [z-1+g, \sigma_{z-1})\cap g\Excl_{g+\ell}(\tau)=[z-1+g, \tau_z)\cap g\Excl_{g+\ell}(\tau)\neq \emptyset. 
  $$
  
  In both cases, we have  deduced that $[z+g-1, \tau_z) \cap g\Excl_{g+\ell}(\tau)\neq \emptyset$,  completing the proof of   Claim 1. 
  
  Claim 1 ensures that we have $\tau\in \mathfrak{S}_n-\mathfrak{S}_{n, g, \ell}$.   
  Next, we proceed to show that the map  $\delta_{n, g,\ell}  $ is an injection.
  By (\ref{deltare2}), we have $c=g\den_\ell(\tau)-g\den_\ell(\pi)$. Therefore,  in order to  show that $\delta_{n, g, \ell}$ is injective, it suffices to show that we can recover $\pi$ from $\tau$.  To this end, we  give  a description of the procedure  to recover the permutation $\pi$ from $\tau$.
  \begin{itemize}
  	 \item Choose $x'$ to be the greatest integer such that the index $x'$ is a $g$-gap $(g+\ell)$-level grande-fixed place of $\tau$.
  	 \item If $n-(g-1)\leq x'\leq n$, let $\pi$ be the permutation obtained from $\tau$ by removing $n$ from $\tau$.
  	   \item If $x'\leq n-g$,   we  generate a permutation  $\sigma$  from $\tau$ by   removing $\tau_{x'}$ from $\tau$. Find all the $g$-gap $(g+\ell)$-level excedance letter of $\sigma$ that are greater than $\tau_{x'}$, say $f_1, f_2, \ldots, f_m$ with $\tau_{x'}<f_1<f_2<\cdots<f_m$. Define $\pi$ to be the permutation obtained from $\sigma$ by replacing $f_i$ with $f_{i-1}$ for all $1\leq i\leq m$ with the convention that $f_0=\tau_{x'}$. 
  \end{itemize} 
 Claim 1 guarantees    that after applying  the above procedure to $\tau$, we  will recover the  permutation $\sigma$ and the permutation permutation $\pi$. This completes the proof. \qed

    \noindent{\bf Proof of Theorem \ref{thden}}.  
   Let $\pi\in \mathfrak{S}_{n-1}$  and $0\leq c\leq n-1$.
   Assume that the space labeled by $c$ in the $g$-gap-$\ell$-level-den-labeling of the index sequence of $\pi$ is immediately before the index $x$ when $c\neq 0$. Then define
   $$  
   \phi^{\den}_{n,g,\ell}(\pi,c)=\left\{ \begin{array}{ll}    
   	\alpha_{n,g,\ell}(\pi, c)&\, \mathrm{if}\,\,   x\in \mathcal{A}(\pi),\\
   	\beta_{n,g,\ell}(\pi, c)&\, \mathrm{if}\,\,  x\in \mathcal{B}(\pi),\\
   	\gamma_{n,g,\ell}(\pi, c)&\, \mathrm{if}\,\,  x\in \mathcal{C}(\pi),\\
   	\delta_{n,g, \ell}(\pi, c)&\, \mathrm{otherwise}. 
   \end{array}
   \right.
   $$
   By the rules specified in the $g$-gap-$\ell$-level-den-labeling of the index sequence of $\pi$, it is easily seen that $c\leq  g\exc_{\ell}(\pi)+g+\ell-2$ when $x\in \mathcal{B}(\pi)$ or $x\in \mathcal{C}(\pi)$, and $c> g\exc_{\ell}(\pi)+g+\ell-2$ when $x\in \mathcal{A}(\pi)$. 
  Then,  by Lemmas \ref{lemalpha1}, \ref{lembeta1},  \ref{lemgamma1} and  \ref{lemdelta1},  we can see that the map   $\phi^{\den}_{n,g,\ell}$ maps the  pair $(\pi, c)$ to a permutation $\tau\in \mathfrak{S}_{n}$ such that the resulting permutation $\tau$ verifies   (\ref{eqprop1}) and (\ref{eqprop2}).  
   By cardinality reasons,
   in order to prove that $\phi^{\den}_{n,g,\ell}$ is a bijection,
   it is sufficient to prove that $\phi^{\den}_{n,g,\ell}$ is injective.  This can be justified by Lemmas  \ref{lemalpha2},    \ref{lembeta2}, \ref{lemgamma3},  \ref{lemgamma4} and 
    \ref{lemdelta2}, completing the proof. \qed

 Let $a=a_1a_2\ldots a_n$ be a weakly increasing sequence of  integers   with $a_i\geq i$ for all $1\leq i\leq n$. Given a permutation $\pi=\pi_1\pi_2\ldots \pi_n\in \mathfrak{S}_{n}$, let $\Exc^{a}(\pi)$ be the subword of $\pi$ consisting of letters $\pi_i$ with $\pi_i>a_i$, in the order induced by $\pi$, and  let $\Nexc^{a}(\pi)$ be the subword of $\pi$ consisting of letters $\pi_i$ with $\pi_i\leq a_i$, in the order induced by $\pi$. Define 
 $$
 \den^{a}(\pi)=\sum\{i+B^a_i(\pi)\mid \pi_i>a_i\}+\inv(\Exc^a(\pi))+\inv(\Nexc^a(\pi)),
 $$ 
  where $B_i^a(\pi)=|\{j\mid 1\leq j<\pi_i\leq a_{j}\}|$.
Han \cite{Han2} proved that $\den^a$  is Mahonian. As remarked by Liu \cite{Liujcta},
 for $n, g, \ell\geq 1$,   we have $\den^a(\pi)=g\den_{\ell}(\pi)$ when $a=a_1a_2\ldots a_n$  with $a_i=\mathrm{max}\{i+g-1, \ell-1\}$ for all $1\leq i\leq n$. This implies that the statistic $g\den_{\ell}$ is Mahonian.

\noindent{\bf Proof of Conjecture \ref{con2}.}  Note that when $n<g+\ell$ and $r=g+\ell-1$, we have
$r\des(\pi)=0=g\exc_{\ell}(\pi)$. From Lemma \ref{Raw}, it is easily seen that the  statistic $r\maj$ is
Mahonian. Since the statistic $g\den_{\ell}$ is also Mahonian, we deduce that   $(r\des, r\maj)$ and $(g\exc_{\ell}, g\den_{\ell})$ are equidistributed over $\mathfrak{S}_{n}$ when $n<g+\ell$. For $n\geq g+\ell$, combining   Lemma \ref{Raw} and Theorem \ref{thden}, we are led to a bijective   proof of Conjecture~\ref{con2}, completing the proof. \qed

\section{Proof of Theorem \ref{thm2}}
This section is devoted to the proof of Theorem  \ref{thm2}. To this end, we need to introduce   a new labeling of the index sequence of a  permutation.

Let $n\geq g+\ell$ and $g,\ell\geq 1$. Given
  $\pi\in \mathfrak{S}_{n-1}$, assume that  $g\Excl_{g+\ell}(\pi)=\{e_1, e_2, \ldots, e_s\}$ with $e_1<e_2<\cdots<e_s$.
  Recall that 
$$
\mathcal{A}(\pi)=\{i\mid    \pi_i\in g\Excl_{g+\ell}(\pi), \,\,     i<\ell\}=\{i\mid \pi_i\geq g+\ell, \,\, i<\ell\}.
$$
Define
$$
\mathcal{B^*}(\pi)=\{i\mid  \pi_i\notin g\Excl_{g+\ell}(\pi), \,\,    i<\ell\}=\{i\mid \pi_i< g+\ell, \,\, i<\ell\}.
$$ 
Recall that $$g\exc_\ell(\pi)=|\{i\mid \pi_i\geq i+g, i\geq \ell\}|=|\{i\mid  \pi_i\in g\Excl_{g+\ell}(\pi), \,\, i\geq \ell\}|.$$  This implies that  
 $|\mathcal{A}(\pi)|+g\exc_\ell(\pi)=s$.  Therefore, we deduce that $|\mathcal{B^*}(\pi)|+s=g\exc_\ell(\pi)+\ell-1$ as $|\mathcal{A}(\pi)|+|\mathcal{B^*}(\pi)|=\ell-1$ and $\mathcal{A}(\pi)\cap\mathcal{B^*}(\pi)=\emptyset $.   
 The {\em $g$-gap-$(g+\ell)$-level-den-labeling} of    the index sequence  of $\pi$ is obtained by labeling the spaces of the index sequence $12\ldots (n-1)$ of $\pi$ as follows.
\begin{itemize}
	\item Label the  rightmost $g$ spaces  by $0,1,\ldots, g-1$ from right to left.
	\item   For all $1\leq i\leq s$,  label the space  before the index $e_i-(g-1) $ with $g+s-i$.
	\item 	Label the spaces  before the indices   belonging to $\mathcal{B^*}(\pi)$ from left to  right with $g+s, \ldots, g\exc_\ell(\pi)+g+\ell-2$.
	\item Label the remaining spaces from left to right with $g\exc_\ell(\pi)+g+\ell-1, \ldots, n-1$.
\end{itemize}
Take    $\pi= 596284317$, $g=3$ and $\ell=6$ for example.
   Clearly, we have $\mathcal{A}(\pi)=\{2\}$,  $\mathcal{B^*}(\pi)=\{1,3, 4,5\}$,    and $3\Excl_{9}(\pi)=\{9\}$. Then the $3$-gap-$9$-level-den-labeling and the $3$-gap $9$-level insertion letters   of the index sequence of  $\pi$ are  given by
   $$
   \begin{array}{lllllllllll} 
   	&_{4}1&_{8}2&_{5}3&_{6}4&_{7}5&_{9}6&_{3}7&_{2}8&_{1}{\bold 9}&_{0}\\
   	&  &  &  &  & &\uparrow &\uparrow&\uparrow&\uparrow&\uparrow\\
   	&  &  &  &  & &9&9&10&10&10 
   	   \end{array}
   $$  
where the labels of the spaces   are written as subscripts and the insertion letters of the spaces  are written on the bottom row.

\begin{remark}\label{re1}
For $x\in \mathcal{A}(\pi)$ or $x\geq \ell$,    the label of the space before the index $x$ under the $g$-gap-$\ell$-level-den-labeling coincides with  the label of the space before the index $x$ under the $g$-gap-$(g+\ell)$-level-den-labeling. 
\end{remark}

Combining  Lemmas \ref{lemalpha1} and \ref{lemalpha2} and Remark \ref{re1}, we derive the following desired properties for $\alpha_{n,g,\ell}$.

\begin{lemma}\label{lemalpha*}
	Let $g, \ell\geq 1$ and $n\geq g+\ell$. The map $\alpha_{n,g,\ell}$ induces an injection from $\mathcal{X}_{n,g, \ell}$ to $\mathfrak{S}^{(1)}_{n, g, \ell}$ such that  for any
	$(\pi, c)\in \mathcal{X}_{n,g, \ell} $,   we have
	\begin{equation}\label{alpha*re1}
		g\exc_\ell(\alpha_{n,g, \ell}(\pi, c))= g\exc_\ell(\pi)+1		 	 
	\end{equation}
	and 
	\begin{equation}\label{alpha*re2}
		g\den_{g+\ell}(\alpha_{n,g, \ell}(\pi, c))=g\den_{g+\ell}(\pi)+c.
	\end{equation}
\end{lemma}

Let $\mathcal{Y}^*_{n,g, \ell}$  denote  the subset of pairs $(\pi, c)\in \mathfrak{S}_{n-1}\times \{ 1, \ldots,  n-1\}$ such that the space labeled by $c$  under the $g$-gap-$(g+\ell)$-level-den-labeling of the index sequence of $\pi$    is immediately  before the index $x$   for some $ x\in\mathcal{B^*}(\pi)$. 
\begin{framed}
	\begin{center}
		{\bf The map $\beta^*_{n,g, \ell}:\mathcal{Y}^*_{n,g, \ell}\longrightarrow  \mathfrak{S}_{n, g, \ell}-\mathfrak{S}^{(1)}_{n, g, \ell} $ }
	\end{center}
	Let $(\pi, c)\in \mathcal{Y}^*_{n,g, \ell}$.  Suppose   that 
	$$  g\Excl_{g+\ell}(\pi)=\{e_1, e_2, \ldots, e_s\}$$ with $e_1<e_2<\cdots<e_s$.  Assume that   the space labeled by $c$  under the $g$-gap-$(g+\ell)$-level-den-labeling of the index sequence of $\pi$    is  immediately before the index $x$.   
	Suppose that $\mathcal{B^*}(\pi)=\{i_1, i_2, \ldots, i_k\}$ with $i_1<i_2<\cdots<i_k$ and that $x=i_{y}$ for some $1\leq y\leq k$. 
	Define  $\beta^*_{n,g,\ell}(\pi,c)$ to be the permutation  constructed by the following procedure.
	\begin{itemize}
		\item Step 1: Construct a permutation  $\sigma=\sigma_1\sigma_2\ldots \sigma_{n-1}$  from $\pi$ by replacing $e_j$ with  $e_{j+1}$ for all $1\leq j\leq s$ with the convention that $e_{s+1}=n$. 
		\item Step 2: Construct a permutation  
		$\theta=\theta_1\theta_2\ldots \theta_{n}$   from $\sigma$ by    replacing $\sigma_x=\sigma_{i_y}$ with $e_1$,
		replacing $\sigma_{i_j}$ with $\sigma_{i_{j-1}}$ for all $y<j\leq k$,   and inserting $ \sigma_{i_k}$ immediately before $\sigma_{\ell}$.	
		
		\item  Step 3: Find all the $g$-gap $(g+\ell)$-level excedance letters of $\theta$ that occur to the left of $\theta_\ell$, say $e'_1, e'_2, \ldots, e'_m$ with  $e'_1<e'_2<\cdots< e'_m$.  
		Then define   $ \beta^*_{n, g, \ell}(\pi,c)$ to be the permutation obtained from $\theta$ by replacing
		$e'_j$ with $e'_{j-1}$ for all $1< j\leq m$ and replacing $e'_1$ with $e'_m$. 
	\end{itemize}
	
\end{framed}
Take    $\pi= 596284317$, $g=3$ and $\ell=6$ for example.
Clearly, we have  $\mathcal{B^*}(\pi)=\{1,3, 4,5\}$   and $3\Excl_{9}(\pi)=\{9\}$. Recall that the $3$-gap-$9$-level-den-labeling    of the index sequence of  $\pi$ is given by 
$$
\begin{array}{lllllllllll} 
	&_{4}1&_{8}2&_{5}3&_{6}4&_{7}5&_{9}6&_{3}7&_{2}8&_{1}{\bold 9}&_{0}.\\
\end{array}
$$
 It is easily seen that the space labeled by $4$ in the   $3$-gap-$9$-level-den-labeling of the index sequence of $\pi$ is immediately before the index $x=1$ and $x\in \mathcal{B}^{*}(\pi)$.  First we generate a permutation $\sigma=5(10)6284317$ from $\pi$ by replacing $e_1=9$ with $e_2=10$.   Then we generate a permutation $\theta= 9(10)56284317 $ from $\sigma$ by replacing $\sigma_1=5$ with $e_1=9$,  replacing $\sigma_3=6$ with $\sigma_1=5$, replacing $\sigma_4=2$ with $\sigma_3=6$,  replacing $\sigma_5=8$ with $\sigma_4=2$,  and inserting $\sigma_5=8$ immediately before $\sigma_6=4$.  
Clearly, all the $3$-gap $9$-level excedance letters that occur to the left of $\theta_6=4$ in $\theta$ are given by $e'_1=9$ and $ e'_2=10$.
Finally, we get the permutation  $\beta^*_{10, 3,6}(\pi,4)=  (10)956284317   $   obtained from $\theta$ by   replacing $e'_2=10$ with $e'_1=9$ and   replacing $e'_1=9$ with $e'_2=10$.

By the same reasoning   as in the proof of Lemmas \ref{lembeta1} and \ref{lembeta2},  one can obtain  following analogous  properties for  $\beta^*_{n,g, \ell}$.  
\begin{lemma}\label{lembeta*}
	Let $g, \ell\geq 1$ and $n\geq g+\ell$. The map  $\beta^*_{n,g, \ell}$ is  an injection from $\mathcal{Y}^*_{n,g, \ell}$ to $\mathfrak{S}_{n, g, \ell}-\mathfrak{S}^{(1)}_{n, g, \ell}$ such that  for any
	$(\pi, c)\in \mathcal{Y}^*_{n,g, \ell} $,   we have
	\begin{equation}\label{beta*re1}
		g\exc_\ell(\beta^*_{n,g, \ell}(\pi, c))= g\exc_\ell(\pi) 	 	 
	\end{equation}
	and 
	\begin{equation}\label{beta*re2}
		g\den_{g+\ell}(\beta^*_{n,g, \ell}(\pi, c))=g\den_{g+\ell}(\pi)+c.
	\end{equation}
\end{lemma}

 Combining Lemmas \ref{lemdelta1} and \ref{lemdelta2} and Remark \ref{re1}, we deduce the   following property of $\delta_{n,g,\ell}$.
  \begin{lemma}\label{lemdelta*}
	Let $g, \ell\geq 1$ and $n\geq g+\ell$.  The map $\delta_{n,g,\ell}$ induces an injection from $\mathcal{W}_{n,g, \ell}$ to $\mathfrak{S}_{n}-\mathfrak{S}_{n,g,\ell}$ such that  for any
	$(\pi, c)\in \mathcal{W}_{n,g, \ell} $, we have  
	\begin{equation}\label{delta*re1}
		g\exc_\ell(\delta_{n,g, \ell}(\pi, c)) =\left\{ \begin{array}{ll}
			g\exc_{\ell}(\pi)&\, \mathrm{if}\,\, 0\leq c\leq g\exc_{\ell}(\pi)+\ell+g-2,\\
			g\exc_{\ell}(\pi)+1&\, \mathrm{otherwise}, 
		\end{array}
		\right.
	\end{equation}
	and 
	\begin{equation}\label{delta*re2}
		g\den_{g+\ell}(\delta_{n,g, \ell}(\pi, c))=g\den_{g+\ell}(\pi)+c.
	\end{equation}
\end{lemma}

\begin{theorem}\label{thden2}
	Fix  $g, \ell\geq 1$ and $n\geq g+\ell$.
	There is a bijection
	$$\psi^{\den}_{n,g,\ell}:  \mathfrak{S}_{n-1}\times \{0, 1, \ldots, n-1\} \longrightarrow \mathfrak{S}_{n}$$ such that  for $\pi\in \mathfrak{S}_{n-1}$ and $0\leq c\leq n-1$,  we have 
	
	\begin{equation}\label{eqprop*1}
		g\exc_\ell(\psi^{\den}_{n,g,\ell}( \pi, c))=\left\{ \begin{array}{ll}
			g\exc_\ell(\pi)&\, \mathrm{if}\,\, 0\leq c\leq g\exc_\ell(\pi)+g+\ell-2,\\
			g\exc_\ell(\pi)+1&\, \mathrm{otherwise}, 
		\end{array}
		\right.
	\end{equation}
	and 
	\begin{equation}\label{eqprop*2}
		g\den_{g+\ell}(\psi^{\den}_{n,g,\ell}( \pi, c))=g\den_{g+\ell}(\pi)+c.
	\end{equation} 
\end{theorem}
\pf 
Let $\pi\in \mathfrak{S}_{n-1}$  and $0\leq c\leq n-1$.
Assume that the space labeled by $c$ in the $g$-gap-$\ell$-level-den-labeling of the index sequence of $\pi$ is immediately before the index $x$ when $c\neq 0$. Then define
$$  
\psi^{\den}_{n,g,\ell}(\pi,c)=\left\{ \begin{array}{ll}    
	\alpha_{n,g,\ell}(\pi, c)&\, \mathrm{if}\,\,   x\in \mathcal{A}(\pi),\\
	\beta^*_{n,g,\ell}(\pi, c)&\, \mathrm{if}\,\,  x\in \mathcal{B^*}(\pi),\\ 
	\delta_{n,g,\ell}(\pi, c)&\, \mathrm{otherwise}. 
\end{array}
\right.
$$
By the rules specified in the $g$-gap-$(g+\ell)$-level-den-labeling of the index sequence of $\pi$, it is easily seen that $c\leq  g\exc_{\ell}(\pi)+g+\ell-2$ when $x\in \mathcal{B}^*(\pi)$, and $c> g\exc_{\ell}(\pi)+g+\ell-2$ when $x\in \mathcal{A}(\pi)$. 
Then,  by Lemmas \ref{lemalpha*}, \ref{lembeta*} and  \ref{lemdelta*},  we can see that the map   $\psi^{\den}_{n,g,\ell}$ is an injection which maps the  pair $(\pi, c)$ to a permutation $\tau\in \mathfrak{S}_{n}$ such that the resulting permutation $\tau$ verifies   (\ref{eqprop*1}) and (\ref{eqprop*2}).  
By cardinality reasons, the map  $\psi^{\den}_{n,g,\ell}$ is indeed a bijection, completing the proof. \qed

\noindent{\bf Proof of Theorem \ref{thm2}.}
 Note that when $n<g+\ell$ and $r=g+\ell-1$, we have
$r\des(\pi)=0=g\exc_{\ell}(\pi)$. Recall that the statistic $r\maj$ is
Mahonian. It is apparent  that  the statistic $g\den_{g+\ell}$ reduces to the statistic $\inv$ when $n<g+\ell$. Therefore,  we deduce that   $(r\des, r\maj)$ and $(g\exc_{\ell}, g\den_{g+\ell})$ are equidistributed over $\mathfrak{S}_{n}$ when $n<g+\ell$.  For $n\geq g+\ell$, combining   Lemma \ref{Raw} and Theorem \ref{thden2}, we are led to a bijective   proof of Theorem ~\ref{thm2}. \qed

\section*{Acknowledgments}
  The authors thank   Zhicong Lin  for his enlightening discussions, from which this work has benefited a lot. 
The work  was supported by
the National Natural
Science Foundation of China grants 12471318 and 12071440.

\end{document}